\documentclass[preprint,preprintnumbers,amsmath,amssymb,floatfix]{revtex4}

\usepackage[dvips]{graphicx}
\usepackage{amssymb}
\usepackage{amsmath}

\newcommand{\ha}{\frac{1}{2}}
\newcommand{\nn}{\nonumber}
\newcommand{\qed}{\hfill \mbox{\raggedright \rule{.07in}{.1in}}}

\def\bO{\mathbf(0)}
\def\bO{\boldsymbol{0}}
\def\eps{\varepsilon}
\def\sov{\omega}

\def\bv{\bar{v}}
\def\bw{\bar{w}}
\def\bn{\bar{n}}
\def\bbv{\bar{\bar{v}}}
\def\bbw{\bar{\bar{w}}}
\def\bbn{\bar{\bar{n}}}
\def\bbbv{\bar{\bar{\bar{v}}}}
\def\bbbw{\bar{\bar{\bar{w}}}}
\def\bbbn{\bar{\bar{\bar{n}}}}
\def\Bf{\bar{f}}
\def\Bg1{\bar{g_1}}
\def\bg2{\bar{g_2}}
\def\bbf{\bar{\bar{f}}}
\def\Bbg1{\bar{\bar{g_1}}}
\def\bbg2{\bar{\bar{g_2}}}
\def\bbbf{\bar{\bar{\bar{f}}}}
\def\Bbbg1{\bar{\bar{\bar{g_1}}}}
\def\bbbg2{\bar{\bar{\bar{g_2}}}}

\newtheorem{thm}{Theorem}
\newtheorem{prop}{Proposition}
\newtheorem{conj}{Conjecture}
\newtheorem{rem}{Remark}

\begin{document}

\preprint{}

\title{A novel canard-based mechanism for mixed-mode oscillations in a neuronal model}

\author{Jozsi Jalics}
\email{jalics@math.ysu.edu}
\affiliation{Department of Mathematics and Statistics, Youngstown State
University, Youngstown, OH, 44555}%

\author{Martin Krupa}%
\email{mkrupa@math.nmsu.edu}
\affiliation{Department of Mathematical Sciences, New Mexico State University, Las Cruces, NM, 88003}%

\author{Horacio G. Rotstein}
\email{horacio@oak.njit.edu}
\affiliation{Department of Mathematical Sciences, New Jersey Institute of Technology, Newark, NJ, 07102}%


\date{\today}

\begin{abstract}
We analyze a biophysical model of a neuron from the entorhinal cortex that includes persistent sodium and slow potassium as non-standard currents using reduction of dimension and dynamical systems techniques to determine the mechanisms for the generation of mixed-mode oscillations. We have found that the standard spiking currents (sodium and potassium) play a critical role in the analysis of the interspike interval.  To study the mixed-mode oscillations, the six dimensional model has been reduced to a three dimensional model for the subthreshold regime.  Additional transformations and a truncation have led to a simplified model system with three timescales that retains many properties of the original equations, and we employ this system to elucidate the underlying structure and explain a novel mechanism for the generation of mixed-mode oscillations based on the canard phenomenon.  In particular, we prove the existence of a special solution, a singular primary canard, that serves as a transition between mixed-mode oscillations and spiking in the singular limit by employing appropriate rescalings, center manifold reductions, and energy arguments.  Additionally, we conjecture that the singular canard solution is the limit of a family of canards and provide numerical evidence for the conjecture. 
\end{abstract}


\maketitle 

\textbf{In this paper we investigate a specific biophysical model of a neuron, not previously analyzed, using techniques from dynamical systems theory.  We discover the presence of  mixed-mode oscillations (MMOs), alternating sequences of large amplitude (neuronal spikes) and small amplitude (subthreshold) oscillations, in the model.  To explain the dynamical mechanisms responsible for the MMOs, we apply reduction of dimensions techniques to reduce the original model to a simplified model system with three timescales.  Our analysis of this system leads to a novel mechanism for the generation of MMOs that is based on the presence of canard solutions that act as a transition between MMOs and spiking.  The theory that we develop may be applied to systems from various contexts that exhibit mixed-mode oscillations and have the same underlying dynamical structure.   
}

\section{Introduction} \label{sec-intro}

Mixed-mode oscillations (MMOs) refer to oscillatory behavior in which 
a number (L) of large amplitude oscillations is followed by 
a number (s) of small amplitude oscillations, 
differing in an order of magnitude, 
to form a complex pattern 
($L_1^{s_1}L_2^{s_2}\dots L_n^{s_n}$) of oscillations of varying amplitudes. 
MMOs were first noticed in the Belousov-Zhabotinsky reaction 
\cite{Zhabotinsky64}, and have been subsequently witnessed in 
experiments and models of numerous chemical 
\cite{Koper95,kn:milszm2,kn:epssho1,Rotstein06b} and biological \cite{Chay95,DicksonOscAct00,Yoshida07}  systems.
Recently there has been work on MMOs in models of
neurons \cite{medvedev04, Drover04, wechselbergerrubin07, Rotstein06,rotsteinwechselberger07, popovic07b}. 
In the context of neural dynamics, large oscillations correspond to firing of action potentials while small oscillations correspond to subthreshold oscillations (STO).

We have discovered MMOs in a biophysical model of a neuron 
in the entorhinal cortex (EC) introduced by Acker et al.\cite{Acker03}.  
The EC is thought to play a prominent role in memory formation through its reciprocal interactions with the hippocampus and participation in the generation of the theta rhythm (4-12 Hz) \cite{Alonso87,Mitchell80,kn:win1,kn:buz2}. 
Pyramidal neurons in layer V of the entorhinal cortex have been found to possess electrophysiological characteristics very similar to those seen in layer II stellate cells of the entorhinal cortex \cite{Schmitz98,Agrawal01}, including their propensity to produce intrinsic subthreshold oscillations at theta  frequencies.  In layer V pyramidal cells these subthreshold oscillations are thought to arise from an interplay between a persistent sodium current and a slow potassium current \cite{Yoshida07,DicksonOscAct00}, while in layer II stellate cells,  they are thought to arise from an interplay between a persistent sodium current and a hyperpolarization activated inward current ($I_h$) \cite{kn:dicalo3,kn:alokli2,Rotstein06}.

The model from Acker et al.\cite{Acker03} that we employ includes these currents and captures some dynamical aspects of the generation of STOs in layer V pyramidal cells. 
After a numerical investigation of the dynamics of the full model that displays the interaction of the persistent sodium and slow potassium currents to form STOs, we modify it by replacing the persistent 
sodium current by a tonic applied current. 
This can be interpreted as the modeling of
an electrophysiological experiment designed to investigate
the effects of blocking the persistent sodium current.  
The modified model displays MMOs, and the main focus of this work is to
understand their nature. In order to simplify the analysis of the model, 
we perform a reduction of dimensions analysis \cite{Rotstein06,kn:clekop1} to reduce the 6-d system 
into a 3-d system. 

Among the proposed mechanisms for MMOs are break-up of an invariant torus
\cite{LS}, break up/loss of stability of a Shilnikov homoclinic orbit
\cite{arne,Koper95} and subcritical Hopf bifurcation combined with a suitable
return mechanism \cite{Guckenheimer1, Guckenheimer2}. 
Well known in the context systems 
with multiple time scales is the {\em canard mechanism}, 
introduced (to the best of our
knowledge) by Milik and Szmolyan \cite{kn:milszm2}. It is based
on the presence of {\em canards}, which are defined as solutions that originate
in a stable slow manifold and continue to an unstable slow manifold while
passing through a thin region of non-hyperbolic behavior, 
which we refer to as the {\em fold region}.
Primary canard solutions, that neither perform small amplitude oscillations 
nor spike, can act as dynamic boundaries between regions of small 
oscillations and spiking.
Canard solutions have been analyzed in three dimensions  
using non-standard analysis by Benoit \cite{Benoit90}
and using techniques of geometric singular perturbation theory in combination 
with blow-up techniques 
by Szmolyan and Wechselberger \cite{Szmolyan01,Wechselberger05}.
A general mechanism of MMOs based on the existence of canards has been 
postulated in \cite{BKW06}.   
There is a variety of different examples which, with some modifications,
fit into the context of this general mechanism 
\cite{Wechselberger05,Drover04,kn:milszm2,
medvedev04, wechselbergerrubin07, Rotstein06, 
rotsteinwechselberger07, popovic07a, popovic07b}.

Our model fits into the context of a system with multiple timescales; 
we will show that the occurrence of MMOs in our model is related to the 
presence of canards. However, while in
most of the other examples of canard mechanism
small oscillations are due to rotations about a canard solution called the
weak canard \cite{BKW06}, in our case, during the stage
of small oscillations, some of the trajectories follow closely 
the unstable manifold of an equilibrium of oscillatory type,
very close to a Hopf bifurcation. In
this aspect our system is similar to the systems studied in \cite{Koper95,
Guckenheimer1, Guckenheimer2}. Consequently our model provides a novel
mechanism for MMOs, one in which both canards and dynamics near a 
saddle-focus equilibrium play an important role.

We designed a three time scale model system 
that accurately reproduces most of the features of the neuronal model, and 
a significant part of this work is devoted to the analysis of the
three time scale model equation. Most of the work, at least to some degree,
relies on numerics. We show that canards can exist,
but we conjecture they occur on parameter intervals of exponentially
small width ($O(e^{-c/\eps})$, where $\eps$ is the ratio of the
two slowest timescales). This is similar to the case of two dimensional
canard explosion (known also as the canard phenomenon), where a small
oscillation grows through a sequence of canard cycles to a relaxation
oscillation as the control parameter moves through an interval of
exponentially small width \cite{Krupa01, Dumortier96, Eckhaus83}.  
In contrast, in the two slow/one fast variable systems studied
by Szmolyan, Wechselberger and others, canards occur in a robust manner,
i.e., on open parameter intervals of size independent of the singular
parameter. We also discover the role of canards in our system:
we show (relying to some degree on numerics) that the transition 
from MMOs to pure spiking must correspond
to a canard solution. We compute a canard numerically and
compare the relevant parameter value to the parameter value of the
transition from MMOs to spiking, finding that the two values are almost equal.
We also calculate analytically the $\eps=0$ approximation of the
numerically computed canard.

An outline of the paper is the following. In section \ref{sec-model}, 
we introduce the model and display the dynamical properties and 
bifurcation structure of the original Acker et al.
model through simulations of the full system.   We also perform 
a reduction of dimensions analysis to 
reduce the 6-d system to a 3-d system with three time scales in the subthreshold regime.
In section \ref{MMOsection}, we introduce the modification of replacing the 
persistent sodium current with a tonic applied current, 
and we analyze the mechanism for the MMOs observed in the modified model. 
We begin with a phase space description of the MMO activity 
through simulations.
Next we propose a three time scale model problem, by adding an 
additional slower time scale 
to the model of Szmolyan and Wechselberger \cite{Szmolyan01}, 
and transform our system into this form.
The three time scale model has a Hopf bifurcation and some
of the dynamics in MMO sequences follows closely
the unstable manifold of the unstable equilibrium,
which reproduces well the behavior of the neuronal model.
In our analysis of the MMOs in the three time scale model, 
a fundamental observation is
that during the approach to the fold region the trajectories are
attracted to a single trajectory, a 1-d (super) slow manifold contained 
in the attracting branch of a  folded 2-d slow manifold.
To investigate the flow near the fold curve, we perform an $\eps$ (singular perturbation parameter) dependent rescaling which blows up the flow near the fold curve. 
We define {\em singular canards} as solutions of 
this system satisfying certain asymptotic properties
and find explicitly a singular canard that is a polynomial 
vector function of the independent variable. 
By employing appropriate rescalings, center manifold reductions, 
and energy arguments, we prove that a continuous family of canard 
solutions must converge, as $\eps\to 0$,
to a singular canard. Consequently, we prove that if a transition from MMO to
spiking occurs along a continuous curve $\mu(\eps)$ then 
$\mu(0)$ must be a value
for which a singular canard exists. We conjecture that the mentioned 
polynomial solution is the only existing singular canard and provide numerical 
evidence for this conjecture.
We conclude with a discussion of related and future work in section 
\ref{sec-discussion}.

\section{The model}\label{sec-model}
\subsection{Model equations}
\begin{figure}[!htb] 
\begin{center}
\includegraphics[angle=270, width=.7\textwidth]{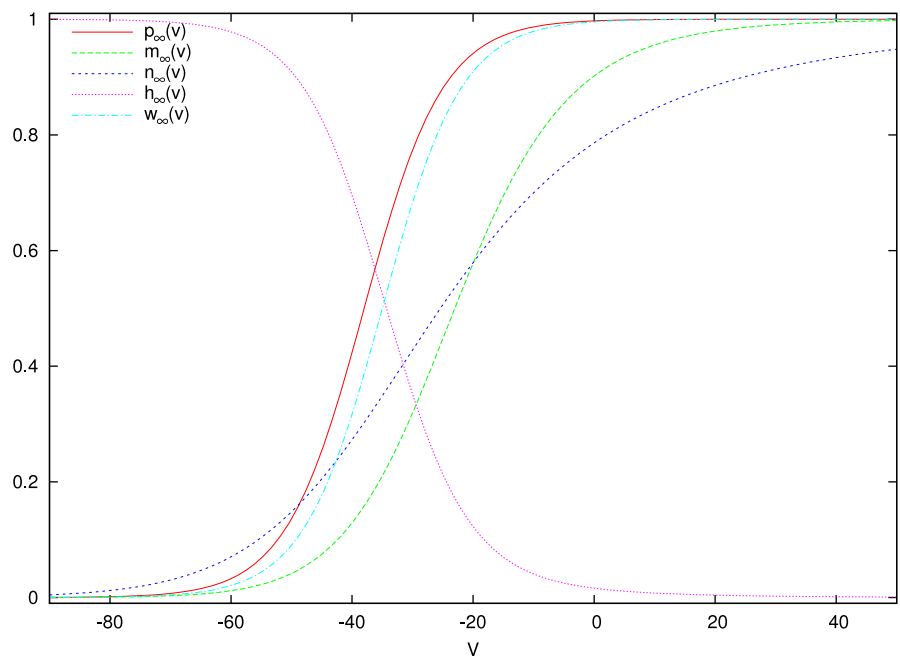}
\caption{Steady state activation and inactivation curves for the gating variables.}\label{figactivCurves}
\end{center} 
\end{figure}

To investigate STOs and MMOs in EC layer V pyramidal cells, 
we consider a biophysical model introduced by Acker et al. \cite{Acker03} that
includes two currents (persistent sodium ($I_{Nap}$) and slow, non-inactivating potassium ($I_{Ks}$) ) whose
interaction is known to be responsible for the generation of STOs in these cells \cite{Yoshida07,DicksonOscAct00}. 
In addition, the model contains the standard Hodgkin-Huxley currents \cite{Hodgkin52}: transient sodium ($I_{Na}$), delayed rectifier potassium ($I_K$), and leak ($I_L$).  The current balance equation for the model is given by 

\begin{eqnarray} \label{fullsystemV}
C\frac{dv}{dt}&=&-g_{Na}m^3h(v-E_{Na})-g_K n^4(v-E_K)- g_L(v-E_L)\\
 && -g_{Nap}p(v-E_{Na})-g_{Ks}w(v-E_K)+ I_{app}  \nn
\end{eqnarray}
where C is the membrane capacitance ($\mu F/cm^2$), $v$ is the membrane potential (mV), $I_{app}$ is the applied bias current ($\mu A/cm^2$), $g$ is conductance ($mS/cm^2$), $E$ is the reversal potential (mV), and the time $t$ is in $ms$.  
Also, we define $I_{Na}=g_{Na}m^3h(v-E_{Na})$, $I_K=g_K n^4(v-E_K)$, $I_L=g_L(v-E_L)$, $I_{Nap}=g_{Nap}p(v-E_{Na})$,
$I_{Ks}=g_{Ks}w(v-E_K)$.   The values of the parameters can be found in the Appendix.   The dynamics of each of the gating variables $m,h,n,p,$ and $w$ is described by a first order equation of the form
\begin{equation}\label{fullsystemGates}
\frac{dx}{dt}= \frac{x_{\infty}(v)-x}{\tau_x(v)}
\end{equation}
 where $x$ represents the gating variable, $\tau_x(v)$ represents the voltage-dependent time constant, and $x_{\infty}(v)$ represents the voltage dependent steady state values of the gating variable. We will see that the $I_{Nap}$ and $I_{Ks}$ play an important role in the generation of subthreshold oscillations. 
The activation and inactivation curves $x_{\infty}(v)$ for the gating variables are shown in Figure \ref{figactivCurves} while the voltage dependent timescales $\tau_x(v)$ are shown in Figure \ref{figRatecurves}.

\subsection{Simulations of the full system}
\begin{figure}[!htb] 
\begin{center}
\includegraphics[angle=270, width=.7\textwidth]{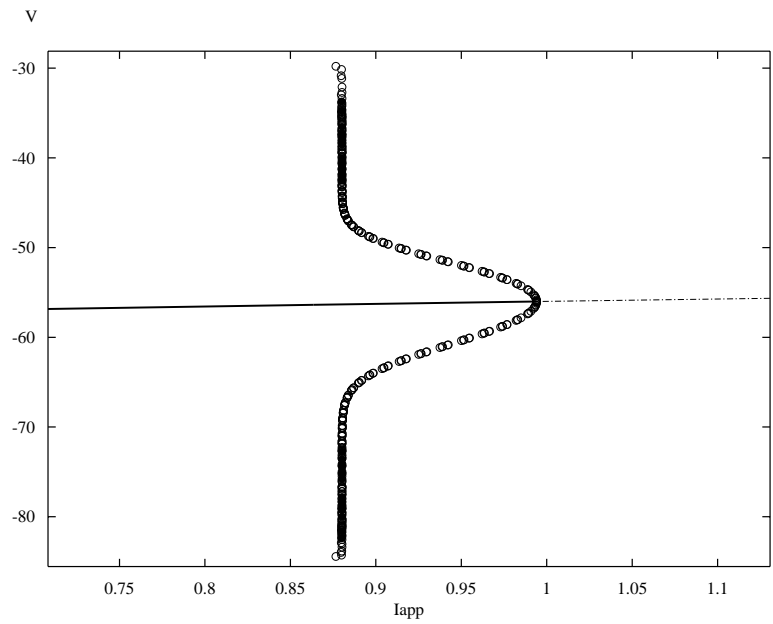}
\caption{Bifurcation diagram for the full system \eqref{fullsystemV}-\eqref{fullsystemGates} as function of $I_{app}$.}\label{figfullsysDefaultBif}
\end{center} 
\end{figure}

\begin{figure}[!htb] 
\begin{center}
\includegraphics[angle=270, width=.7\textwidth]{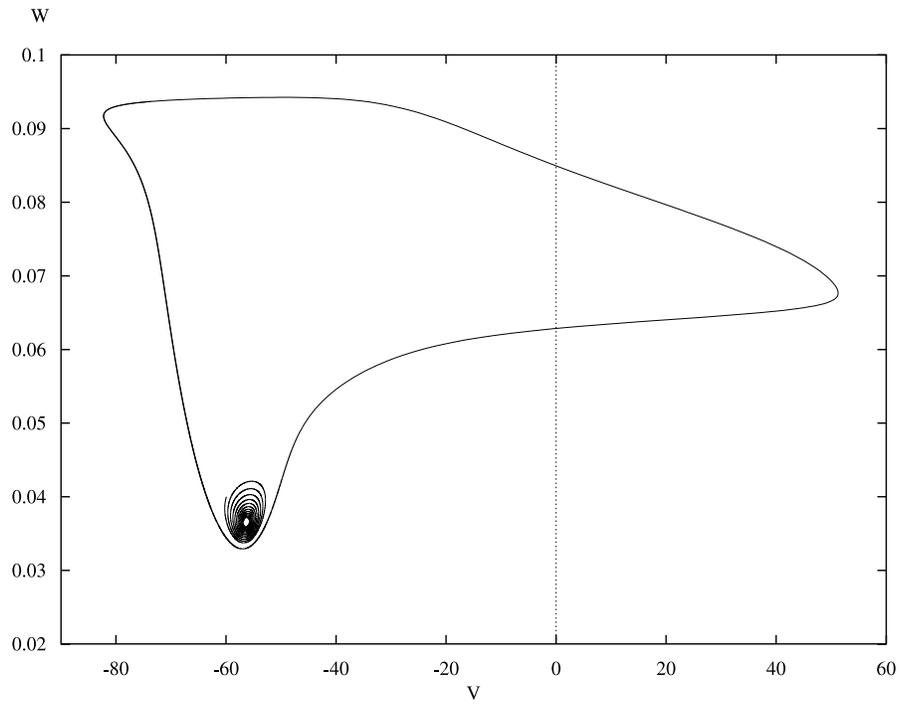}
\caption{Solution trajectories of the full system \eqref{fullsystemV}-\eqref{fullsystemGates} with $I_{app}=.9$ projected onto the $w$-$v$ plane showing bi-stability of spiking solutions and damped subthreshold oscillations.}\label{fullbistab}
\end{center} 
\end{figure}

\begin{figure}[!htb] 
\begin{center}
\includegraphics[angle=270, width=.7\textwidth]{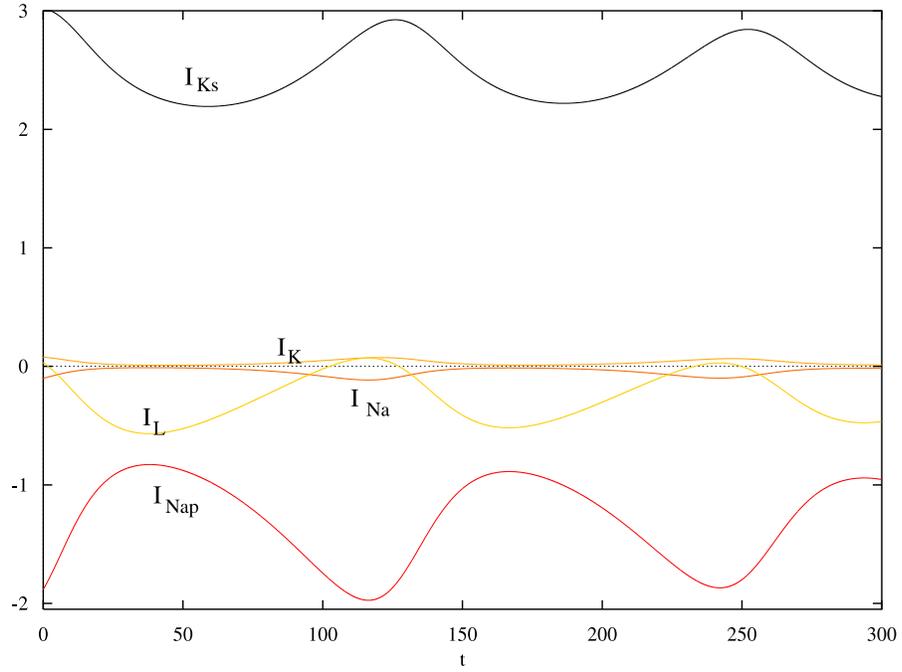}
\caption{Currents in the full system \eqref{fullsystemV}-\eqref{fullsystemGates} with $I_{app}=.9$ during damped subthreshold oscillations showing the interaction of $I_{Nap}$ and $I_{Ks}$.}\label{STOfull}
\end{center} 
\end{figure}

Simulations of the full system \eqref{fullsystemV}-\eqref{fullsystemGates} show that the model fits into the category of class 2 neural excitability \cite{kn:izh2}. Using the applied current as a bifurcation parameter, the system undergoes a subcritical Hopf bifurcation as $I_{app}$ is increased as seen by Figure \ref{figfullsysDefaultBif}. Thus, as $I_{app}$ is increased, the cell switches between three different regimes .  In the first regime ($I_{app}<0.88$), there is a single stable rest state which is converged to through damped subthreshold oscillations.  In the second regime ($0.88 <I_{app}<0.994$),  there is bi-stability of the rest state and periodically spiking solutions, which are separated by unstable, small amplitude oscillations generated through a subcritical Hopf bifurcation (See Figure \ref{fullbistab}). The spiking solutions start with a minimal frequency of about 5Hz while the damped subthreshold oscillations have a frequency of about 10Hz.  Also, the replacement of the slow potassium current with a tonic drive results in the absence of subthreshold oscillations (as in \cite{Yoshida07}) while the replacement of the persistent sodium current with a tonic drive will be discussed in section \ref{MMOsection}.  We note that the damped subthreshold oscillations involve the interplay of the slow potassium current and the persistent sodium current as seen in Figure \ref{STOfull}.  At the transition between the first and second regimes, we observe dynamics reminiscent of a 2-d canard phenomenon 
in which the trajectory shown follows a repelling unstable manifold. 
The third regime ($I_{app}>0.994$) consists of stable, periodic firing.

 \subsection{Reduction of Dimensions}

\begin{figure}[!htb] 
\begin{center}
\includegraphics[angle=270, width=.7\textwidth]{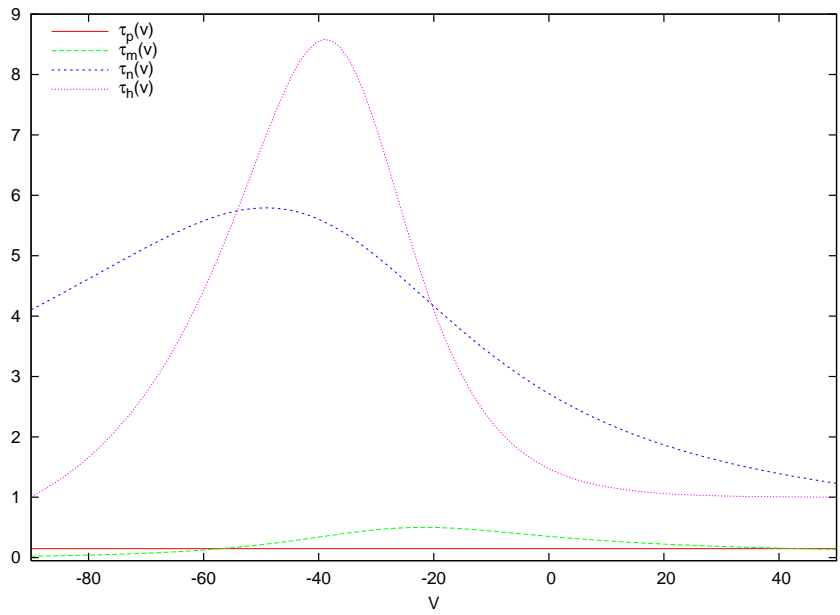}
\caption{Voltage dependent time scale curves for the gating variables. Note: $\tau_w(v)=90$.}\label{figRatecurves}
\end{center} 
\end{figure}

We show that during the interspike interval, which corresponds to the subthreshold voltage regime, we can reduce the system to a 3-dimensional system (\ref{reduced3dsys}).  
By analyzing the voltage dependent timescales $\tau_x(v)$ seen in Figure \ref{figRatecurves}, we see that $\tau_w(v) >> \tau_n(v),\tau_h(v)>> \tau_m(v),\tau_p(v)$. Thus, we seem to have three timescales with $p$ and $m$ evolving very quickly and $w$ evolving very slowly.  Thus, letting $p=p_{\infty}(v)$ and $m=m_{\infty}(v)$ yields a reasonable approximation.  
By considering Figure \ref{figactivCurves} and \ref{figRatecurves}, we see that when $v$ is hyperpolarized, $h$ evolves much faster to its steady state value (its time scale is ~4 times smaller) than $n$. 
Also, we note that $n_{\infty}'(v) > h_{\infty}'(v)$. 
Thus, $h$ quickly approaches its steady state value $h_{\infty}(v)$, and then it can track it closely since $h_{\infty}(v)$ changes very slowly during the interspike interval. Thus, we replace $h$ with $h_{\infty}(v)$ and retain $n$.   Hence, we can reduce the system to one for which only $v,w,$ and $n$ are dynamic variables.   

\begin{figure}[!htb] 
\begin{center}
\includegraphics[angle=270, width=.7\textwidth]{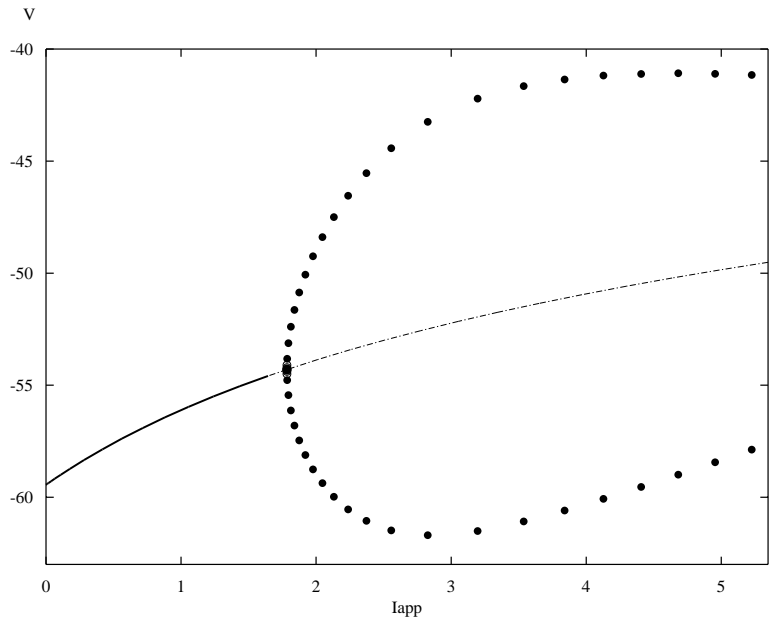}
\caption{Bifurcation diagram for the full system \eqref{fullsystemV}-\eqref{fullsystemGates} with $I_{Na}$ removed as function of $I_{app}$. Note the change in the location and criticality of the Hopf bifurcation point versus Figure \ref{figfullsysDefaultBif}. }\label{figfullsysDefaultgNa0Bif}
\end{center} 
\end{figure}

\begin{figure}[!htb] 
\begin{center}
\includegraphics[angle=270, width=.7\textwidth]{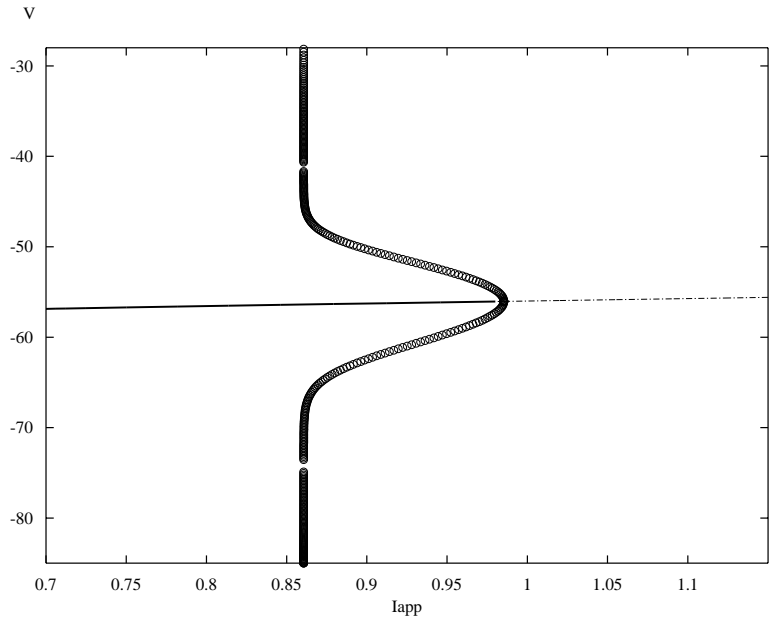}
\caption{Bifurcation diagram for the reduced 3-d system \eqref{reduced3dsys} as a function of $I_{app}$.}\label{3dNbif}
\end{center} 
\end{figure}

To investigate the importance of the spiking currents (sodium and potassium), we also perform a bifurcation analysis.
When $I_{Na}$ is removed, the system exhibits a supercritical Hopf bifurcation in which sustained, stable subthreshold oscillations appear.  (See Figure \ref{figfullsysDefaultgNa0Bif}.)  This change in the bifurcation criticality due to $I_{Na}$ removal shows that the transient sodium current plays an important role in the interspike interval.
With $I_{K}$ removed, the bifurcation structure is retained, but there is a significant shift in the Hopf bifurcation point (from $I_{app}\approx .9941$ to $I_{app} \approx .775$). 
By retaining the $n$ equation, and replacing $m,h,$ and  $p$ with their steady state activation and inactivation functions $m_{\infty}(v),h_{\infty}(v),$ and $p_{\infty}(v)$, respectively, we obtain the reduced 3d system

\begin{eqnarray}\label{reduced3dsys}
C\frac{dv}{dt}&=&-g_{Na}m_{\infty}^3(v)h_{\infty}(v)(v-E_{Na})-g_K n^4(v-E_K)- g_L(v-E_L) \nn\\
&& -g_{Nap}p_{\infty}(v)(v-E_{Na})-g_{Ks}w(v-E_K)+ I_{app}  \\
\frac{dw}{dt}&=& \frac{w_{\infty}(v)-w}{\tau_w}\nn\\
\frac{dn}{dt}&=& \frac{n_{\infty}(v)-n}{\tau_n(v)}\nn\
\end{eqnarray}
where $\tau_w(v)$ has been replaced by its constant value, $\tau_w$.
We note that the bifurcation diagram for this system (see Figure \ref{3dNbif}) is very similar to that of the full system.  

Typically, the spiking currents are not included in the biophysical description of the subthreshold regime as in 
generalized integrate-and-fire models, which are constructed by biophysically describing the subthreshold regime using non-spiking currents and adding artificial spikes \cite{Richardson03,Rotstein06}.   
This system provides an example where neglecting the spiking currents in the biophysical description of the subthreshold regime leads to a qualitative change in the cell's dynamics, i.e. $I_K$ and $I_{Na}$ must be included in this reduced system since our model is sensitive to small changes in the sodium and potassium currents during the latter part of the interspike interval before the spike. 
In some cases it has been shown that these currents play no critical role in the subthreshold regime \cite{Rotstein06}.

\section{Analysis of mixed-mode oscillations}\label{MMOsection}  

\subsection{Simulations and outline of the results}

\begin{figure}[!htb] 
\begin{center}
\includegraphics[angle=270, width=.7\textwidth]{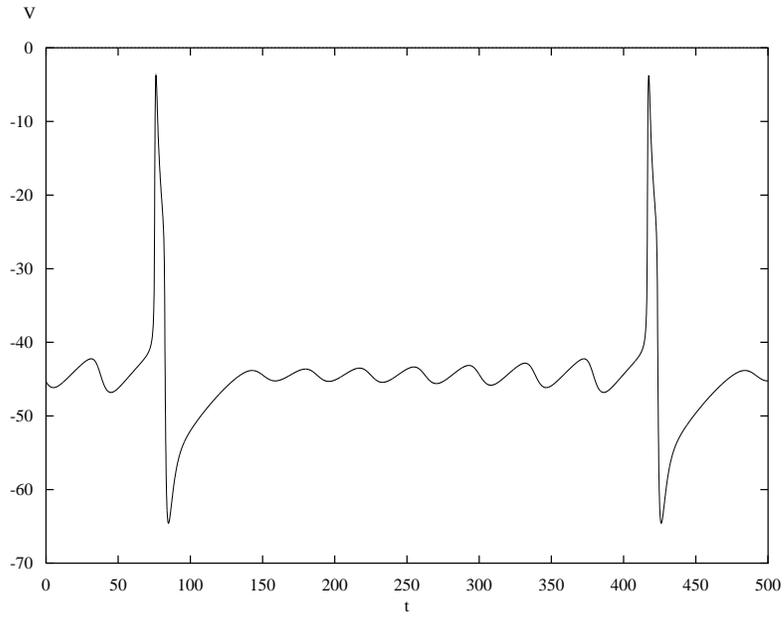}
\caption{Mixed-mode oscillations in the 3-d system \eqref{reduced3dsys} with $I_{app}=17.1$.}\label{deepEKs3dMMOIapp17_1.ode}
\end{center} 
\end{figure}

\begin{figure}[!htb] 
\begin{center}
\includegraphics[angle=270, width=.7\textwidth]{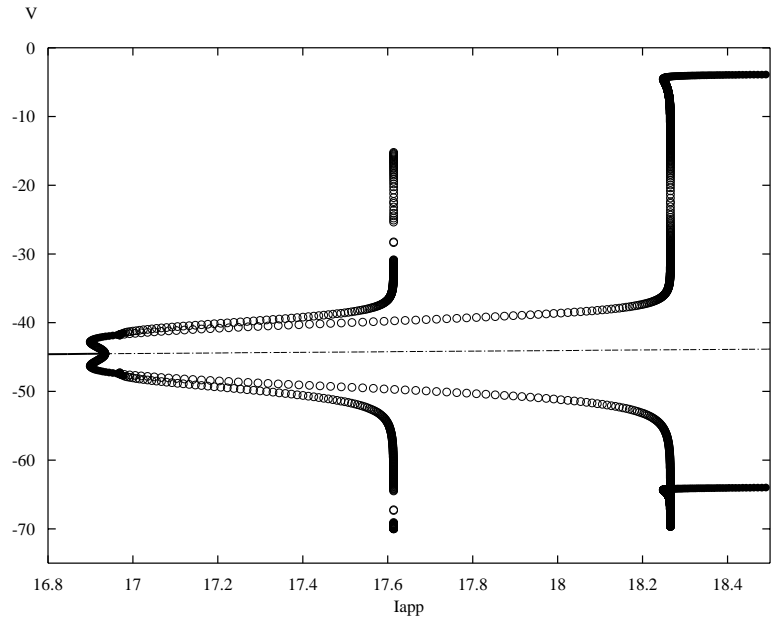}
\caption{Bifurcation diagram for the 3-d system \eqref{reduced3dsys} with the persistent sodium current removed. 
A sequence of period doubling bifurcations occur starting around $I_{app}=16.97$.}\label{3dBif}
\end{center} 
\end{figure}

We have discovered that in the full system (\ref{fullsystemV}-\ref{fullsystemGates}) and in the reduced system (\ref{reduced3dsys}) the replacement of the persistent sodium current with a tonic excitatory drive leads to mixed-mode oscillations in which a large amplitude spike is followed by a number of small amplitude subthreshold oscillations (see Figure \ref{deepEKs3dMMOIapp17_1.ode}).  We note that this modified version of  \eqref{reduced3dsys}, with the persistent sodium current removed, exhibits a more complex bifurcation structure than \eqref{reduced3dsys} itself as seen in Figure \ref{3dBif}. 
The system exhibits a single equilibrium which loses stability through a subcritical Hopf bifurcation at $I_{app} \approx16.93$.   For $I_{app} \gtrapprox18.25$, the system exhibits continuous spiking. For $I_{app}$ between $\approx17.61$ and $\approx18.25$, the system displays MMOs of type $1^1$. For $I_{app}$ between $\approx 16.98$ and $\approx 17.61$, the system displays MMOs with multiple STOs per spike, and the number of STOs per spike is a decreasing function of $I_{app}$. 
Finally, for $I_{app}$ between $\approx16.9$ and $\approx16.98$, the system exhibits sustained subthreshold oscillations.

\begin{figure}[!htb] 
\begin{center}
\includegraphics[width=\textwidth]{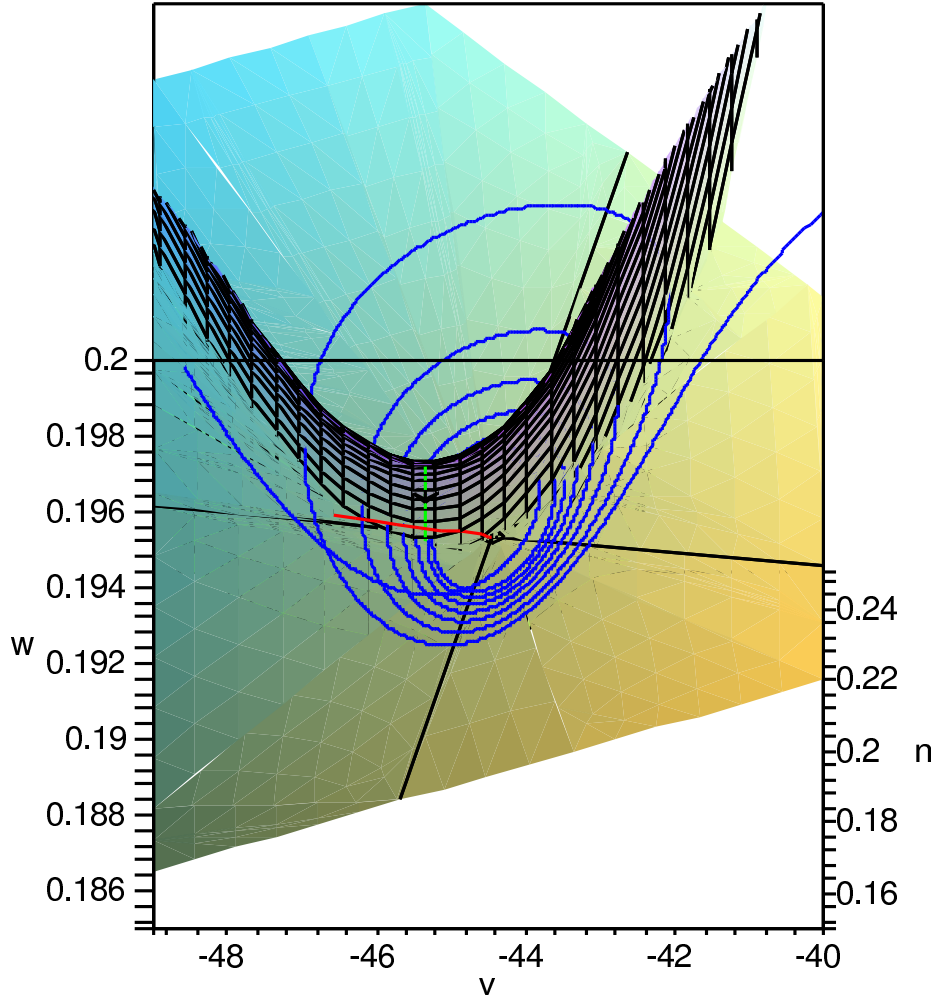}
\caption{MMO trajectory for Iapp=17.1 in \eqref{reduced3dsys}.  Trajectory in blue. Fold curve and folded singularity in green.  Fixed point and its 1-d stable eigenvector in red. Folded surface is the v-nullsurface, and the shaded plane is the unstable eigenspace of the fixed point.
(Note that the trajectory pierces the v-nullsurface near its fold. This is typical for passage near a folded node or folded saddle node type of singularity.)}\label{figorig4}
\end{center} 
\end{figure}

Our goal is to analyze the mechanisms for this behavior in \eqref{reduced3dsys} with $g_{NA}=0$.
We begin with a geometric interpretation of the MMOs using simulations of \eqref{reduced3dsys}.
We note that the  v-nullsurface is a 2-d surface $w(v,n)$, and for fixed n, $w(v,n)$ is a cubic shaped curve.  Also, we note that \eqref{reduced3dsys} possesses a single fixed point, denoted by $P_0=(v_0,w_0,n_0)$, which is near the fold curve of the v-nullsurface.  
The linearization about $P_0$ yields a large negative eigenvalue and a pair of complex eigenvalues with positive real part, and the ratio of the real eigenvalue to the real part of the complex eigenvalue is about -30. Thus,  $P_0$ possesses a 1-d stable manifold and a 2-d unstable manifold.  
Numerical simulations show that a MMO
trajectory travels down near the left sheet of the v-nullsurface (slow manifold)  
and subsequently passes near the fold, where it 
gets attracted to the 2-d unstable manifold 
of the fixed point $P_0$. See Figure \ref{figorig4}.  
The 1-d stable eigenspace of $P_0$ is plotted in red in the figures, and the plane displayed in the figures is the unstable 
eigenspace of $P_0$.   The trajectory approaches the unstable manifold of $P_0$, which we denote by $M_0$, along the stable 
manifold of $P_0$.  Then the trajectory spirals out from near $P_0$ along $M_0$ and continues to transversally pierce the 
v-nullsurface (the loss of alignment with the fast nullsurface during the passage through the fold region
is typical for folded node or folded saddle node type dynamics), until it becomes tangent to the v-nullsurface. 
At this point, the trajectory leaves $M_0$ and gets attracted to the right sheet of the v-nullsurface since it enters a region in which $v$ 
is strongly repelling from the left sheet. Then, there is a return mechanism (spiking regime) that brings the trajectories back to the left sheet of 
the v-nullsurface. We note that the feature of a passage near 
the unstable manifold of a saddle focus equilibrium  followed by a reinjection
mechanism is reminiscent of Shilnikov's example \cite{Shilnikov65} as well as 
of the examples where MMOs occur through a subcritical Hopf bifurcation
\cite{Guckenheimer1, Guckenheimer2} 
(see section \ref{sec-discussion} for a detailed 
discussion).

In the remaining part of this work much of the focus will be on the
following model system
\begin{align}\label{eq-3tssimplesl}
\begin{split}
\eps^2&\dot x=-y+x^2+D_1x^3\\
&\dot y=x-z\\
\eps &\dot z=\mu+ax+bz, 
\end{split}
\end{align}
where $\eps>0$ is small and $\mu$ is a control parameter.
System  \eqref{eq-3tssimplesl} shares many of the features
of \eqref{reduced3dsys} but is a lot simpler. 
The initial part of our analysis is to establish two basic
geometric features of the flow:
\begin{enumerate}
\item
There exists an equilibrium of oscillatory type, which loses stability
through a Hopf bifurcation for $\mu$ close to $0$.
\item
The critical manifold $y=x^2+D_1x^3$ is a surface with two fold lines.
There exists a 2D attracting slow manifold close to one of the
sheets of this  surface, and within this slow manifold there is a 
1D super-slow manifold, approximating the curve $y=x^2+D_1x^3$, $x=z$.
\end{enumerate}
We will show that all the MMO trajectories must pass very close 
to the super-slow manifold.
Consequently,  the distance between
the continuation of the super-slow manifold to the fold region and the
stable manifold of the equilibrium is of great importance for 
the nature of the MMOs. If this distance is small,
then the trajectory must follow for a long time
the flow along the 2D unstable manifold of the equilibrium, making
many rotations about the equilibrium. 
As this distance increases, the role of the unstable manifold of the
equilibrium diminishes.

In this work we have focused on the transition from MMOs to spiking.
This occurs in the parameter regime where the continuation of the super-slow
manifold, and thus the MMO trajectories, are relatively far from the equilibrium. 
We have been able to show, with the help of some numerics, 
that the transition from MMOs to spiking occurs
through a canard, and that families of canard solutions converge
to singular canards, which are special solutions of a suitable
limiting system. We have also found a parameter value for which
a special solution exists. We conjecture, based on numerical evidence,
that this parameter value is unique.
\subsection{Derivation of the model system}
In order to investigate MMOs in system \eqref{reduced3dsys}, 
we transform it to a three timescale toy model of the form
\begin{eqnarray}\label{canon3d}
x'&=& -y+ F(x)\nn\\
y'&=& \eps^2(x-z) \\
z'&=& \eps(\mu+ax+bz) \nn\
\end{eqnarray}
with $F(x)=x^2+D_1x^3$ and $\varepsilon\ll 1$. 
To this end, we begin by first treating \eqref{reduced3dsys} as a slow-fast system with two timescales in order to obtain a folded singularity, which we will use as the origin for our model system in order to study small amplitude oscillations near the fold curve of the 2-d folded v-nullsurface of \eqref{reduced3dsys}.  For an exposition on folded singularities, see  \cite{Szmolyan01,Wechselberger05}.

Two timescales are uncovered by dividing the $v$ equation of (\ref{reduced3dsys}) by the maximum conductance $g_{Na}$ and setting $\kappa=C/g_{Na}$ we obtain 
\begin{eqnarray}\label{reduced3dsys2}
\kappa \frac{dv}{dt}&=&-m_{\infty}^3(v)h_{\infty}(v)(v-E_{Na})-\bar{g}_K n^4(v-E_K)- \bar{g}_L(v-E_L) \nn\\
&& -\bar{g}_{Nap}p_{\infty}(v)(v-E_{Na})-\bar{g}_{Ks}w(v-E_K)+ \bar{I}_{app}\equiv f(v,w,n)  \\
\frac{dw}{dt}&=& \frac{w_{\infty}(v)-w}{\tau_w} \equiv g_1(v,w,n) \nn\\
\frac{dn}{dt}&=& \frac{n_{\infty}(v)-n}{\tau_n(v)} \equiv g_2(v,w,n) \nn\
\end{eqnarray}
where $\bar{g}_x=g_x/g_{Na}$ for $x=K,L,Nap$, and $Ks$ and $\bar{I}_{app}={I}_{app}/g_{Na}$. We note that the conductances in the $v$ equation are bounded by 1 and the magnitudes of the $w$ and $n$ equations are bounded by 1 since $\tau_w=90$ and $1\le\tau_n(v)\le 6$.    
Thus, $v$ evolves on a faster time scale then $w$ or $n$.  Momentarily, we neglect the fact that $w$ and $n$ vary on different timescales.

By setting $\kappa=0$, we obtain the reduced system
which is a 2-d dynamical system on the 
critical manifold $S=\{(v,w,n):f(v,w,n)=0\}$.  Based on the work of \cite{Szmolyan01,Wechselberger05}, we wish to show that $S$ is a non-degenerate folded surface in the neighborhood of a folded equilibrium point $(v_c,w_c,n_c)$ with 
$f(v_c,w_c,n_c)=0,f_v(v_c,w_c,n_c)=0,f_w(v_c,w_c,n_c)\ne 0,
f_{vv}(v_c,w_c,n_c)=0$. 
In order to compute $(v_c,w_c,n_c)$, we begin by solving
$f(v,w,n)=0$ for $w$ and obtain $w=W(v,n)$.  Next, we solve 
$f_v(v,W(v,n),n)=0$ for $v$ as a function of $n$ and find that $v$ is a constant, which we denote by
$v_c$, along the fold curve (since $f$ and $f_v$ are linear in $n^4$ and $w$).  Also, we denote $W(v_c, n)$ by $\xi(n)$.
Thus, the fold curve is parametrized by $(v_c, \xi(n), n)$, $n \in I$.
Also, we differentiate $f(v,w,n)=0$ with respect to $t$ to obtain
$f_v \dot v=-f_w g_1-f_n g_2$.  Thus, the reduced 2-d system on the 
critical manifold can be denoted by
\begin{eqnarray}\label{reduced2dsys}
-f_v\dot{v}&=& f_w g_1+f_n g_2 \\
\dot{n}&=& g_2 \nn\
\end{eqnarray}
where $w$ is replaced by $W(v,n)$.
Through a solution dependent time rescaling we
obtain the desingularized reduced system
\begin{eqnarray}\label{reduced2dsys2}
\dot{v}&=&f_wg_1+f_ng_2  \\
\dot{n}&=& -f_vg_2.\nn\
\end{eqnarray}
We numerically solve $f_w g_1+f_n g_2=0$ for $n$ with $v=v_c$ and $w=
\xi(n)$ to obtain $n_c$. 
Finally, $w_c= \xi(n_c)$.  We note that $(v_c,w_c,n_c)$
is a fixed point of (\ref{reduced2dsys}) and (\ref{reduced2dsys2}).
We linearize (\ref{reduced2dsys2}) about this fixed point and find that in our parameter regime both eigenvalues are \emph{negative}, resulting in a folded node equilibrium. We note that the quotient of these eigenvalues predicts the maximum number of subthreshold oscillations per spike \cite{Szmolyan01,Wechselberger05} for $\kappa$ sufficiently small.   However, this prediction is incorrect (underestimate) for this system since there are three timescales, resulting in a $\kappa$ that is not sufficiently small in the two timescale formulation.    
 
We now shift the folded node point to the origin in \eqref{reduced3dsys2}
with the transformation  $\bv=v-v_c$, $\bw=w-w_c$, and $\bn=n-n_c$ to obtain
  \begin{eqnarray}\label{3dtranf1}
\kappa \dot{\bv} &=& f(\bv+v_c,\bw+w_c,\bn+n_c)\equiv\Bf(\bv,\bw,\bn)   \\
\dot{\bw}&=& g_1(\bv+v_c,\bw+w_c,\bn+n_c) \equiv\Bg1(\bv,\bw,\bn) \nn\\
\dot{\bn}&=& g_2(\bv+v_c,\bw+w_c,\bn+n_c) \equiv\bg2(\bv,\bw,\bn). \nn\
\end{eqnarray}  
Next, we rectify the fold curve to the $n$ axis.  We note that the fold curve is parametrized by $(0, \chi(\bn), \bn)$. where $\chi(n)= W(v_c,\bn+n_c)-w_c$. Thus, we make the substitution
$\bbv=\bv$, $\bbw=\bw-\chi(\bn)$, and $\bbn=\bn$ to obtain
 \begin{eqnarray}\label{3dtranf2}
\kappa \dot{\bbv} &=& \Bf(\bbv,\bbw +\chi(\bbn),\bbn)  \equiv\bbf(\bbv,\bbw,\bbn)  \\
\dot{\bbw}&=& \Bg1(\bbv,\bbw +\chi(\bbn),\bbn)-\chi_n(\bbv,\bbw +\chi(\bbn),\bbn)\bg2(\bbv,\bbw +\chi(\bbn),\bbn)   \equiv \Bbg1(\bbv,\bbw,\bbn) \nn\\
\dot{\bbn}&=& \bg2(\bbv,\bbw +\chi(\bbn),\bbn)  \equiv\bbg2(\bbv,\bbw,\bbn). \nn\
\end{eqnarray} 
We note that $\bbf(\bO)=0=\Bbg1(\bO)$.

Then we let $\bbbv=\gamma_1\bbv$, 
$\bbbw=\gamma_2\bbw$, and $\bbbn=\bbn$
where $\gamma_1=\ha\bar{\bar f}_{vv}(\bO)$ and $\gamma_2=\ha\bar{\bar {f}}_{vv}(\bO)\bar{\bar {f}}_{w}(\bO)$ to obtain

. \begin{eqnarray}\label{3dtranf3}
\kappa \dot{\bbbv} &=& \bbbv+\bbbw^2+ O(\bbbv\bbbw,\bbbv^3)\equiv \bbbf(\bbbv,\bbbw,\bbbn)  \\
\dot{\bbbw}&=& \Bbbg1(\bbbv,\bbbw,\bbbn) \nn\\
\dot{\bbbn}&=& \bbbg2(\bbbv,\bbbw,\bbbn). \nn\
\end{eqnarray} 
To simplify the $n$ equation we let $v=\bbbv$, $w=\bbbw$, and $n=\bbbn/\bbbg2(\bO)$ to obtain
\begin{eqnarray}\label{3dtranf4}
\kappa \dot{v}&=& v^2+w +O(vw,v^3)\nn\\
\dot{w}&=& \beta n-cv + O(w,n^2,nv,v^2) \\
\dot{n}&=& 1 + dn+ev +O(w,v^2,nv). \nn\
\end{eqnarray}
Next, we rescale time with $t=\tau/\kappa$, let $x=v, y=-w,$ and $z=(\beta/c)n$, and truncate higher order terms to obtain
\begin{eqnarray}\label{canon3dlinear3}
x'&=& -y+x^2 +D_1x^3\nn\\
y'&=& {\kappa c}(x-z) \\
z'&=& \frac{\kappa\beta}{c}(1+\frac{cd}{\beta}z+ex). \nn\
\end{eqnarray}
Finally, if we let $\eps=\sqrt{\kappa c}$, then we obtain
\begin{eqnarray}\label{canon3dlinear4}
x'&=& -y+x^2 +D_1x^3\nn\\
y'&=& \eps^2(x-z) \\
z'&=& \eps(\mu+ax+bz) \nn\
\end{eqnarray}
  where $\mu=\beta \sqrt{\kappa} c^{-3/2})$, $a=\beta e \sqrt{\kappa}c^{-3/2}$, and $b=d \sqrt{\kappa}c^{-3/2}$.  This system is equivalent to \eqref{canon3d}.  For example, if $I_{app}=17.1$ in \eqref{reduced3dsys}, the parameters are approximated by: $\eps =.2764436178\equiv \eps_d,\; \mu =.02839202004,\; a =.4266759683$, and $b =-.9420074624 $.
Also, we consider $D_1=-.4$, which provides an appropriate return mechanism.
A more general assumption that we make to guarantee the desired functioning
of the toy model is 
\begin{equation}\label{eq-ab}
-b>a>0.
\end{equation}   
\subsection{Slow manifolds of the model system}\label{sec-S1D}
\begin{figure}[!htb] 
\begin{center}
\includegraphics[width=.95\textwidth]{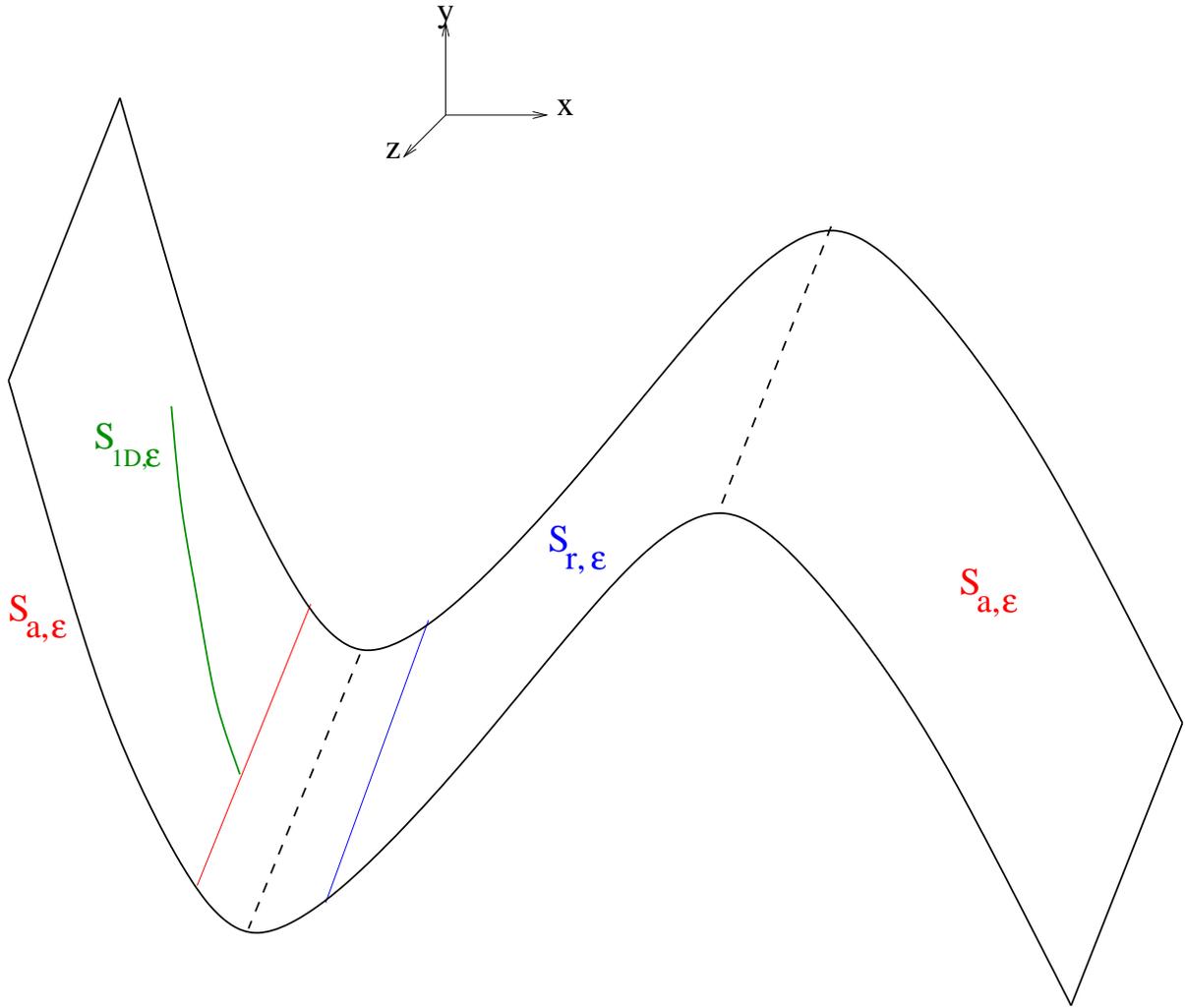}
\caption{Caricature of the folded slow manifold.}\label{foldMMO}
\end{center} 
\end{figure}

Equation \eqref{canon3d} has a two dimensional {\em critical manifold}
defined as the set of equilibria of \eqref{canon3d} with $\eps$ set
to $0$. We denote this manifold by $S$ and note that it is defined
by the equation $y=x^2+D_1x^3$. Let $x_{max}$ denote the unique value
of $x$ for which the function $F(x)=x^2+D_1x^3$ has a maximum.
Let
\begin{align*}
\begin{split}
&S_a=\{(x,y,z)\in S\; :\; x\le 0\}\\
&S_r=\{(x,y,z)\in S\; :\; 0\le x\le x_{max}\}
\end{split}
\end{align*}
where $a$ and $r$ denote attracting and repelling, respectively.
By Fenichel theory \cite{Fenichel71}, away from the fold lines $x=y=0$ and
$x=x_{max}$, $y=f(x_{max})$, system \eqref{canon3d} with $\eps>0$ has
locally invariant manifolds $S_{a,\eps}$ and $S_{r,\eps}$ that are $O(\eps)$ close
to $S_a$ and $S_r$. While the existence of $S_{a,\eps}$ and $S_{r,\eps}$
is standard and not surprising, we also argue that there exists 
a one-dimensional super-slow manifold (it can be chosen to be contained in 
$S_{a,\eps}$), denoted by $S_{1D,\eps}$, that attracts all the dynamics on $S_{a,\eps}$. (See Figure \ref{foldMMO}.)  
For this reason it is convenient to rescale
\eqref{canon3d} to the slow time scale, which gives
\eqref{eq-3tssimplesl}.
Setting $\eps=0$, we obtain two constraints:
\begin{align*}
\begin{split}
0&=-y+x^2+D_1x^3\\
0&=\mu+ax+bz
\end{split}
\end{align*}
This gives a one dimensional curve defined by:
\begin{align}\label{eq-1dsm}
\begin{split}
y&=x^2+D_1x^3\\
z&=-\frac{\mu +ax}{b},
\end{split}
\end{align}
which we denote by $S_{1D}$. Note that assumption \eqref{eq-ab} guarantees that
the flow on $S_{1D}$ is towards the fold.

\begin{figure}[!htb] 
\begin{center}
\includegraphics{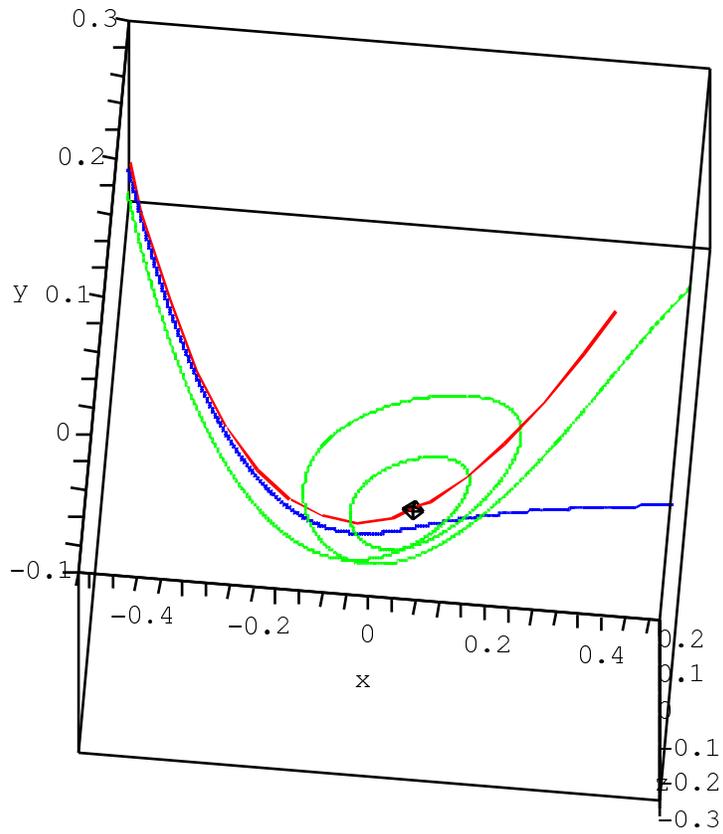}
\caption{ $S_{1D}$ plotted along with a spiking solution of \eqref{canon3d} with $\eps=.1$ that tracks $S_{1D}$ quite closely as well as an MMO solution for the larger value of $\eps=\eps_d=.2764436178$. $\mu=.02839$. }\label{slowMan1d}
\end{center} 
\end{figure}

Linearizing \eqref{eq-3tssimplesl}
about $S_{1D}$ we obtain the matrix
\begin{align*}
A=\left (\begin{array}{ccc}\frac{1}{\eps^2}(2x+3D_1x^2)&-\frac{1}{\eps^2}&0\\
                            1&0&-1\\
                            \frac{a}{\eps}&0&\frac{b}{\eps}\end{array}\right ).
\end{align*}
Ignoring the $O(1)$ terms in $A$ we obtain
\begin{align*}
\left (\begin{array}{ccc}\frac{1}{\eps^2}(2x+3D_1x^2)&-\frac{1}{\eps^2}&0\\
                            0&0&0\\
                            \frac{a}{\eps}&0&\frac{b}{\eps}\end{array}\right ).
\end{align*}
It is easy to verify that this matrix, apart from an eigenvalue $0$, has eigenvalues
$\frac{1}{\eps^2}(2x+3D_1x^2)$ and $\frac{b}{\eps}$, which are 
both negative for $x<0$. 
It follows that the dominant eigenvalues of $A$ are
\[
\frac{1}{\eps^2}(2x+3D_1x^2)+O(1)\qquad\mbox{and}\qquad \frac{b}{\eps}+O(1).
\]
There is also a third eigenvalue, which is $O(1)$, but we are interested
only in the two leading ones.
By a slight modification of Fenichel theory we can
prove that \eqref{canon3d} has a one dimensional attracting locally
invariant manifold, which we denote by $S_{1D,\eps}$. 
All the trajectories which make a relaxation loop
(spike) must become exponentially close to it before they enter 
the fold region.  Thus, there is a single trajectory that needs to be tracked through the fold region to determine the existence 
of MMOs/small amplitude dynamics
and this trajectory can be arbitrarily chosen as long as its initial condition
is not too close to the fold.  
We illustrate this scenario in Figure \ref{slowMan1d} in which we display $S_{1D}$ 
along with a spiking solution for $\eps=.1$ that tracks $S_{1D}$ quite closely as well as a MMO solution for the larger value of $\eps=\eps_d$.

\subsection{The microscope: a rescaling of the model system}\label{sec-micro}
In order to understand the flow of \eqref{canon3d} in the fold region,
we introduce a rescaling which blows up a rectangle of size
$O(\eps)\times O(\eps^2)\times O(\eps)\times O(\eps)$ to the size of $O(1)$ 
in each of the directions. This trick has been widely used 
in the  study of canard problems
and was given the name 'microscope' by Benoit, Diener and others
who investigated canard explosion and related problems by means of
non-standard analysis \cite{Benoit81,Benoit90,Diener94}.
In the case of  (\ref{canon3d}), the 'microscope' rescaling is as follows:
\[
x=\eps\bar x, y=\eps^2\bar y, z=\eps\bar z, \mu=\eps\bar\mu .
\]
\begin{rem}
{\rm In Proposition \ref{prop-bound} we will prove that MMOs
cannot exist if $\bar\mu$ is too large. This implies that MMOs exist
for $\mu=O(\eps)$ only.}
\end{rem}
After rescaling the variables we obtain the
following system: 
\begin{align}\label{eq-rescD1}
\begin{split}
\bar x'=&-\bar y+\bar x^2 + \eps D_1 \bar x^3\\
\bar y'=&\bar x-\bar z\\
\bar z'=&\bar \mu+a\bar x+b\bar z.\\
\end{split}
\end{align}
By setting $\eps = 0$ we obtain:
\begin{align}\label{eq-resc}
\begin{split}
\bar x'=&-\bar y+\bar x^2\\
\bar y'=&\bar x-\bar z\\
\bar z'=&\bar \mu+a\bar x+b\bar z.\\
\end{split}
\end{align}
\begin{rem}\label{rem-fn}
{\rm A similar normal form arises in the study of folded node 
\cite{Benoit90}, \cite{Szmolyan01}. The difference is in
the $\bar z'$ equation, where the terms $a\bar x$ and $b\bar z$
do not appear. However these two terms alter the dynamics dramatically,
so that the results on the folded node normal form cannot be easily
applied in the context of \eqref{eq-resc}. We have some partial results
in the direction of reducing the study of the dynamics of
\eqref{eq-resc} to facts about the dynamics of the normal form of folded
node, but the arguments are very technical and it does not seem worthwhile
to present these results in this paper.}
\end{rem}
The system \eqref{eq-resc} has a unique equilibrium given by
\[
\bar x=\bar z=-\frac{\bar\mu}{a+b},\quad \bar y=\left (\frac{\bar\mu}{a+b}\right )^2.
\]
The linearization at the equilibrium is given as follows:
\[
A_0=\left (\begin{array}{ccc}-\frac{2\bar\mu}{a+b}&-1&0\\
                           1&0&-1\\
                           a&0&b\end{array}\right ) 
\]
with ${\rm det}\, A_0= a+b$.
We assume that $(a+b)<0$, which implies that the determinant is always
negative. Clearly the only way for the matrix to be non-hyperbolic is if 
there is a pure imaginary eigenvalue. Let $p(\lambda)$ denote 
the characteristic  polynomial of $A_0$ and note that $p(\lambda)$ can be
written in the form $p(\lambda)=\lambda^3+a_2\lambda^2+a_1\lambda+a_0$.
The condition for $A_0$ to have a purely imaginary eigenvalue is
$a_0=a_1a_2$. We evaluate the coefficients and get
\begin{align*}
a_0=-(a+b),\quad a_1=\frac{-2\bar\mu b+a+b}{a+b},
\quad a_2=\frac{-b(a+b)+2\bar\mu}{a+b}.
\end{align*}
We define the function $C(\bar\mu)$ by the formula
\[
C(\bar\mu)=(a_0-a_1a_2)(a+b).
\]
Clearly the roots of $C(\bar\mu)$ give purely imaginary eigenvalues.
Evaluating $C(\bar\mu)$ we get
\[
C(\bar\mu)=(a+b)a+(2b^2+2)\bar\mu-\frac{4b}{a+b}\bar\mu^2.
\]
Now $C(\bar\mu)$ has two roots, given by
\begin{align*}
(1+b^2\pm \sqrt{(1+b^2)^2+4ba})\frac{a+b}{4b}.
\end{align*}
For example, if $a=.4266759683$ and $b=-.9420074624$ (corresponding to $I_{app}=17.1$), the two roots are  
$\bar\mu=0.06692624850$ and $\bar\mu=0.4493251582$.  The first corresponds to the Hopf bifurcation point while the second is an extraneous solution.

\subsection{Special solution}\label{sec-spec}
System \eqref{eq-resc}  has special solutions that are algebraic as 
$t\to -\infty$, i.e. solutions that tend to $-\infty$ no faster than a polynomial function of $t$. In this section we will prove that for each value of 
$\bar\mu$ there exists exactly one such solution (up to time rescaling). 
The importance of special solutions will become clear  
in Section \ref{sec-cont} where we will show that the continuation
of the manifold $S_{1D,\eps}$ to the region of validity of the 'microscope coordinates'
$(\bar x,\bar y, \bar z)$ is $O(\eps)$ close to the special solution
of \eqref{eq-resc} corresponding to the same value of $\bar\mu$.

In order to understand the flow of \eqref{eq-resc} near infinity, 
including special solutions, we make the following transformation:
\begin{align}\label{eq-trafo}
\sigma=\frac{1}{\sqrt{\bar y}},\quad\tilde x=\sigma\bar x,\quad \tilde z= \sigma\bar z.
\end{align}
The transformation yields the following system:
\begin{align}\label{eq-rescc}
\begin{split}
&\tilde x'=\sigma^{-1}(-\frac12\sigma^2(\tilde x-\tilde z)\tilde x
-1+\tilde x^2)\\
&\sigma '=-\frac12\sigma^2(\tilde x-\tilde z)\\
&\tilde z'=\sigma\bar\mu+a\tilde x+ b\tilde z
-\frac12\sigma(\tilde x-\tilde z)\tilde z.
\end{split}
\end{align}
Through an appropriate (solution dependent) time rescaling we can eliminate
a factor of $\sigma^{-1}$ in \eqref{eq-rescc}:
\begin{align}\label{eq-resccr}
\begin{split}
&\tilde x'=-\frac12\sigma^2(\tilde x-\tilde z)\tilde x-1+\tilde x^2\\
&\sigma '=-\frac12\sigma^3(\tilde x-\tilde z)\\
&\tilde z'=\sigma(\sigma\bar\mu+a\tilde x+ b\tilde z
-\frac12\sigma(\tilde x-\tilde z)\tilde z).
\end{split}
\end{align}
Solutions of \eqref{eq-resccr}, suitably reparametrized, become solutions of 
\eqref{eq-rescc}, except of course for the solutions for which $\sigma=0$.
Equation \eqref{eq-resccr} has two lines of equilibria,
given by
\begin{align}\label{eq-le}
\tilde x=\pm 1,\quad\sigma=0.
\end{align}
These lines of equilibria were artificially introduced by the time rescaling 
which led from \eqref{eq-rescc} to \eqref{eq-resccr}. The reason for 
making the transformation was to remove the singularity created by the term 
$\sigma^{-1}$. The resulting system \eqref{eq-resccr} can be analyzed 
using dynamical systems theory.  Namely,  we
can understand the flow near the line of equilibria, thus obtaining insight 
into the structure of solutions of \eqref{eq-rescc}. We focus on the line of
equilibria with $\tilde x=-1$ since it corresponds to the attracting $x$ 
direction. Observe that there exists a center manifold containing the line 
of equilibria. We denote this center manifold by $M_{2D}$. Note also that 
points on $M_{2D}$ satisfy the constraint $\tilde x=-1+O(\sigma^2)$
by center manifold theory \cite{Carr81}.  
Substituting into \eqref{eq-resccr}, we obtain
\begin{align}\label{eq-rcrcm}
\begin{split}
&\sigma '=-\frac12\sigma^3(-1-\tilde z+O(\sigma^2))\\
&\tilde z'=\sigma(\sigma\bar\mu-a+ b\tilde z + O(\sigma^2)-\frac12\sigma(-1-\tilde z+O(\sigma^2))\tilde z).
\end{split}
\end{align}
Now we can introduce the "inverse" of the time rescaling which lead 
from \eqref{eq-rescc} to \eqref{eq-resccr}, by reducing a factor of 
$\sigma$. Thus we are back to the original time coordinate. The resulting 
system is
\begin{align}\label{eq-rcrcmrb}
\begin{split}
&\sigma '=-\frac12\sigma^2(-1-\tilde z+O(\sigma^2))\\
&\tilde z'=\sigma\bar\mu-a+ b\tilde z+O(\sigma^2)-\frac12\sigma(-1-\tilde z+O(\sigma^2))\tilde z.
\end{split}
\end{align}
System \eqref{eq-rcrcmrb} has an equilibrium at 
$(\sigma=0,\tilde z=a/b)$, which
has one stable direction and one center direction. The branch of the 
center manifold existing for $\sigma>0$, which we denote by $M_{1D}$,  
is unique. Equation \eqref{eq-rcrcmrb}
restricted to $M_{1D}$ has the form
\begin{align}\label{eq-M1D}
\sigma '=-\frac12\sigma^2(-1-\frac{a}{b}+O(\sigma))
\end{align}
and gives a solution of \eqref{eq-rescc} and thus of 
\eqref{eq-resc} which is algebraic in $t$ backwards in time. All other 
solutions of \eqref{eq-rescc} grow at least exponentially 
as $t\to -\infty$. We will denote the special solution in the 
$(\bar x, \bar y, \bar z)$ coordinates by $\gamma_{\bar\mu}$.

We now look for a special solution of \eqref{eq-resc} which is a polynomial vector function in $t$.
Suppose that $\gamma_{{\rm poly}}(t)=(\bar x(t),\bar y(t),\bar z(t))$ is a polynomial solution and note that it can only exist if $\bar x(t)$ and $\bar z(t)$ are linear 
in $t$ and $\bar y(t)$ is quadratic in $t$. We can write
\begin{align*}
\bar x(t)=\bar x_0+\alpha t,\qquad \mbox{with $\alpha$ some unknown constant.}
\end{align*}
Using the first equation in \eqref{eq-resc}, we obtain
\begin{align*}
\bar y(t)= (\bar x_0 +\alpha t)^2-\alpha.
\end{align*}
Substitution into the second equation in \eqref{eq-resc} results in
\begin{align*}
\bar z(t)= (1-2\alpha)(\bar x_0+\alpha t).
\end{align*}
Finally, the third equation in \eqref{eq-resc} gives
\begin{align*}
(\alpha-2\alpha^2)=\bar \mu +(a\alpha+ b(\alpha-2\alpha^2))t,
\end{align*}
from which we get the conditions
\begin{align*}
\begin{split}
&\alpha-2\alpha^2=\bar \mu\\
&a\alpha+ b(\alpha-2\alpha^2)=0.
\end{split}
\end{align*}
Solving we get
\begin{align*}
\alpha=\frac{a+b}{2b},\qquad\bar\mu=-\frac{a(a+b)}{2b^2}.
\end{align*}
It follows that a polynomial solution exists for $\bar\mu=\bar\mu_0$ with
\[
\bar\mu_0=-\frac{a(a+b)}{2b^2}.
\]
 Note that  transformation \eqref{eq-trafo} takes solutions with algebraic behavior as $t\to-\infty$ to solutions with algebraic behavior as $t\to-\infty$.
Consequently, by uniqueness of
$M_{1D}$, $\gamma_{\bar\mu_0}=\gamma_{{\rm poly}}$. 
We denote the value of $\mu$ corresponding to  $\bar\mu_0$ by $\mu_0$.
\begin{rem} {\rm We could repeat the construction of the beginning 
of this section
to obtain special solutions which are algebraic as $t\to\infty$.
For each $\bar\mu$ there exists exactly one such solution. Typically
the special solutions for $t\to -\infty$ do not match with the
special solutions for $t\to\infty$. The number $\bar\mu_0$ is an exceptional
value of $\bar\mu$ for which the two kinds of special solutions are equal;
the corresponding solution is $\gamma_{{\rm poly}}$.
In Section \ref{sec-percan} we will present numerical evidence
that $\gamma_{{\rm poly}}$ is the $\eps\to 0$ limit of a 
family of canard solutions.}
\end{rem}
\subsection{Continuation of $S_{1D,\eps}$  into the fold region}
\label{sec-cont}
In Section \ref{sec-S1D} we proved that all trajectories 
on $S_{a,\eps}$ are attracted to a one dimensional manifold 
$S_{1D,\eps}$. Naturally $S_{1D,\eps}$ and its flow continuation
play a crucial role in the organization of the dynamics
of \eqref{canon3d}. In this section, we show that any solution
with initial condition on $S_{1D,\eps}$, when it passes through
the region of validity of the $(\bar x, \bar y, \bar z)$ coordinates,
must be $O(\eps)$ close to the special solution
$\gamma_{\bar\mu}$.

As in Section \ref{sec-spec}, we use the coordinate transformation 
\eqref{eq-trafo}, but we apply it to equation \eqref{eq-rescD1} here. This yields
\begin{align}\label{eq-rescD1c}
\begin{split}
&\tilde x'=\sigma^{-1}(-\frac12\sigma^2(\tilde x-\tilde z)\tilde x
+(-1+\tilde x^2+\frac{\eps}{\sigma}D_1\tilde x^3)\\
&\sigma '=-\frac12\sigma^2(\tilde x-\tilde z)\\
&\tilde z'=\sigma\bar\mu+a\tilde x+ b\tilde z
-\frac12\sigma(\tilde x-\tilde z)\tilde z.
\end{split}
\end{align}
Subsequently, we apply the time rescaling resulting in the removal of the factor
$\sigma^{-1}$:
\begin{align}\label{eq-rescD1cr}
\begin{split}
&\tilde x'=-\frac12\sigma^2(\tilde x-\tilde z)\tilde x+(-1+\tilde x^2
+\frac{\eps}{\sigma}D_1\tilde x^3)\\
&\sigma '=-\frac12\sigma^3(\tilde x-\tilde z)\\
&\tilde z'=\sigma(\sigma\bar\mu+a\tilde x+ b\tilde z
-\frac12\sigma(\tilde x-\tilde z)\tilde z).
\end{split}
\end{align}
Notice that setting $\eps=0$ in \eqref{eq-rescD1c} (resp.  \eqref{eq-rescD1cr})
gives \eqref{eq-rescc} (resp. \eqref{eq-resccr}).

We now note that
\[
\sigma=\frac{\eps}{\sqrt{y}},\quad \tilde x=\frac{x}{\sqrt{y}}, 
\quad \tilde z=\frac{z}{\sqrt{y}},
\]
where $(\tilde x,\sigma , \tilde z)$ are the coordinates introduced in
\eqref{eq-trafo}. Consider a point $p\in S_{1D,\eps}$ close to the fold
region and let $\tilde p$ be the corresponding point in the
$(\tilde x,\sigma , \tilde z)$ coordinates. 
By \eqref{eq-1dsm}, $p$ is 
close to the curve $(-\rho,\rho^2, -\frac{\mu-a\rho}{b})$. We fix 
$\rho>0$ small, which corresponds to the assumption 
that $p$ is close to the fold region. Observe that
$\tilde p$ is close to $(-1,\frac{\eps}{\rho} , -\frac{\mu-a\rho}{b\rho})$. 
We investigate the solution of \eqref{eq-rescD1cr} with initial 
condition $\tilde p$.

Consider the following system of equations:
\begin{align}\label{eq-rescD1crp}
\begin{split}
&\tilde x'=-\frac12\sigma^2(\tilde x-\tilde z)\tilde x+(-1+\tilde x^2+D_1\omega\tilde x^3)\\
&\sigma '=-\frac12\sigma^3(\tilde x-\tilde z)\\
&\tilde z'=\sigma(\sigma\bar\mu+a\tilde x+ b\tilde z-\frac12\sigma(\tilde x-\tilde z)\tilde z)\\
&\omega'=\frac12\sigma^2\omega(\tilde x-\tilde z).
\end{split}
\end{align}
Note that \eqref{eq-rescD1crp} has a constant of motion $H(\sigma, \omega)= \omega\sigma$, and that
\eqref{eq-rescD1cr}  is a reduction of \eqref{eq-rescD1crp} to the invariant hypersurface $\eps=\omega\sigma$.

We now argue in a similar way as in section \ref{sec-spec}:
\begin{itemize}
\item there exists a center manifold satisfying the constraint 
$\tilde x= -1+O(\sigma^2, \omega)$. We denote this manifold by $M_{2D}^\eps$ because, 
even though it is three dimensional in the context of \eqref{eq-rescD1crp}, 
it is two dimensional in the variables $(\tilde x, \sigma, \tilde z)$, 
i.e. once the reduction $\omega=\frac{\eps}{\sigma}$ has been implemented.
The manifold $M_{2D}^\eps$ can be viewed as a
center manifold of the 
2D surface of equilibria defined by $\sigma=0$ and $\tilde x= -1+O(\omega)$. 

\item restricted to $M_{2D}^\eps$ we can rescale the time and find a 
2D center manifold, denoted by $M_{1D}^\eps$, which is a natural generalization of $M_{1D}$ to the 
$\eps >0$ case. That is, center manifold theory guarantees the existence of 
functions $\psi_j(\sigma,\omega)$ such that $M_{1D}^\eps$ 
is defined by the equations
\begin{equation}\label{eq-formM1de}
\tilde x=
-1+\psi_1(\sigma, \omega),\qquad \tilde z=\frac{a}{b}+\psi_2(\sigma, \omega),
\end{equation}
and $M_{1D}$ is defined by the equations
\begin{equation}\label{eq-formM1d}
\tilde x=-1+\psi_1(\sigma, 0),\qquad \tilde z=\frac{a}{b}+\psi_2(\sigma, 0).
\end{equation}
(Recall that $M_{1D}$ corresponds to the special solution.)
The flow on $M_{1D}^\eps$ is given by 
\begin{align}\label{eq-rescD1M2D
}
\begin{split}
&\sigma '=-\frac12\sigma^2(-1-\frac{a}{b}+O(\omega,\sigma))\\
&\omega'=\frac12\sigma\omega(-1-\frac{a}{b}+O(\omega,\sigma))
\end{split}
\end{align}
and the level curves of the function $H=\omega\sigma$ are invariant. (This system is the analogue of \eqref{eq-M1D} in the $\eps=0$ case.) 
\end{itemize}

\begin{figure}[!htb] 
\begin{center}
\includegraphics[width=.45\textwidth]{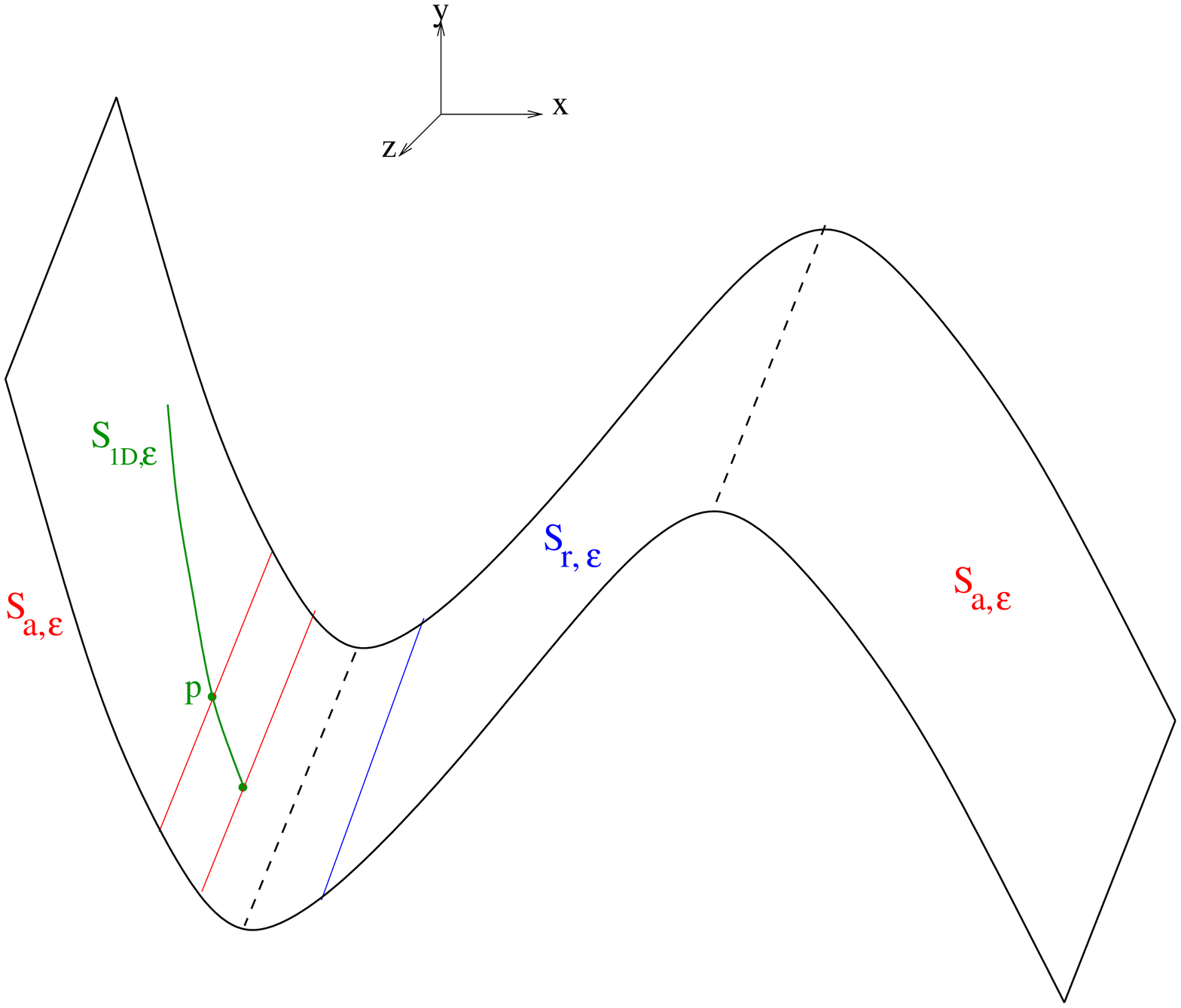}
\includegraphics[width=.45\textwidth]{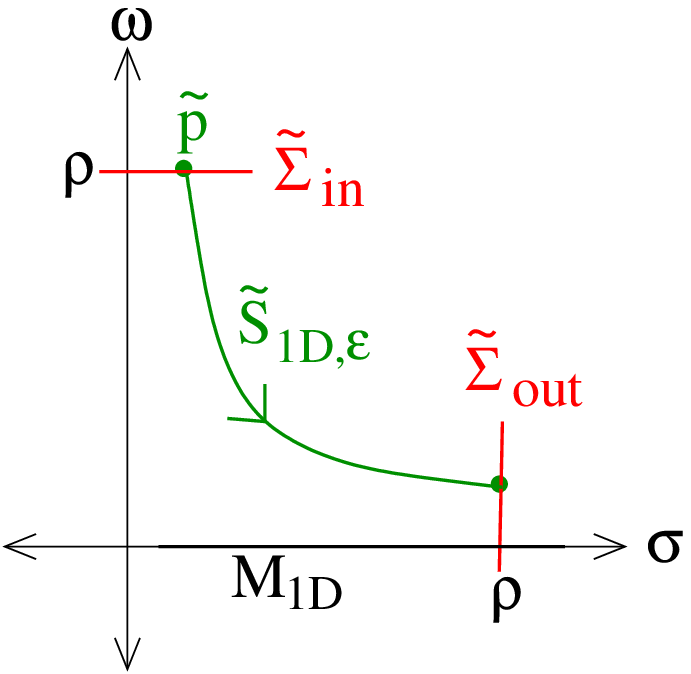}
\caption{Caricature of the flow near the fold curve.}\label{foldMMOsect}
\end{center} 
\end{figure}

We define two sections of the flow of \eqref{eq-rescD1crp}:
\[
\tilde\Sigma^{in}=\{(\tilde x,\sigma ,\tilde z, \rho),\qquad \tilde\Sigma^{out}=\{(\tilde x,\rho ,\tilde z, \omega),
\]
where $\rho$ is the constant used in the definition of the point $p$. (See Figure \ref{foldMMOsect}.)
We are interested in describing the transition
by the flow of \eqref{eq-rescD1crp} from $\tilde\Sigma^{in}$ to $\tilde\Sigma^{out}$.
We state a few important  properties of \eqref{eq-rescD1crp}:
\begin{enumerate}
\item The manifold $M_{2D}^\eps$ is strongly attracting. To see this, recall
that $M_{2D}^\eps$ is an attracting center manifold of the line of equilibria
defined by $\tilde x=-1$, $\sigma=0$.
\item The manifold $M_{1D}^\eps$ is also attracting, since it is
an attracting center manifold in the context of the system obtained by
restricting \eqref{eq-rescD1crp} to $M_{2D}^\eps$ and cancelling a factor
of $\sigma$ (similar construction as leading from \eqref{eq-rcrcm}
to \eqref{eq-rcrcmrb}).
The attraction rate
is weaker, but since the solutions on $M_{1D}^\eps$ move fairly slowly,
solutions   with initial conditions off $M_{1D}^\eps$ in $\tilde\Sigma^{in}$ are very close
to $M_{1D}^\eps$ in $\tilde\Sigma^{out}$. 
\item 
Consider a trajectory $\phi(t)$ of \eqref{eq-rescD1crp} starting in
$\tilde\Sigma^{in}$. Since $\phi(t)$ 
is contained in the invariant hypersurface  
$\omega\sigma =\eps$, it follows that it
must approach the intersection of $\omega\sigma =\eps$ with 
the attracting manifold $M_{1D}^\eps$.
Since this is a one dimensional curve containing no singularities, 
it must be a trajectory. Hence it follows that $\phi(t)$ 
approaches a specific trajectory on $M_{1D}^\eps$.
\item By the definition of $\tilde\Sigma^{out}$ and due to the constraint
$\eps=\sigma\omega$ we have 
$\tilde\Sigma^{out}\cap M_{1D}^\eps =\{
(-1+\psi_1(\rho, \frac{\eps}{\rho}), \rho, 
\frac{a}{b}+\psi_2(\rho, \frac{\eps}{\rho}),\frac{\eps}{\rho})\}$,
which, together with \eqref{eq-formM1d} implies that the distance between
$\tilde\Sigma^{out}\cap M_{1D}^\eps$
and $\tilde\Sigma^{out}\cap M_{1D}$ is $O(\eps)$.
Let $T$ be defined by the requirement $\{\phi(T)\}=
\tilde\Sigma^{out}\cap M_{1D}^\eps$. It follows that
$\phi(T)$ converges to $\gamma_{\bar\mu}$ as $\eps\to 0$.
\end{enumerate}

It follows from points 1. and 2.
that the transition map from  $\tilde\Sigma^{in}$ 
to $\tilde\Sigma^{out}$, defined by the  flow 
of \eqref{eq-rescD1cr}, mapping a neighborhood of $\tilde p$ 
to a neighborhood of $\phi(T)$, is a strong contraction. 
We can conclude this because
the restriction of the flow of \eqref{eq-rescD1crp} 
using the constraint $\eps=\omega\sigma$ to the flow of
\eqref{eq-rescD1cr} amounts to factoring out the neutral direction.
\begin{rem}\label{rem-mr2de}
{\rm
The manifold $M^\eps_{2D}$ as we defined it (center manifold 
of the surface of equilibria given by $\sigma=0$ and $\tilde x= -1+O(\omega)$ for \eqref{eq-rescD1crp}) 
is not unique. We add a requirement that $M^\eps_{2D}$ equals
$S_{a,\eps}$ in the section $\tilde\Sigma^{in}$.
This way $M^\eps_{2D}$ is the continuation of $S_{a,\eps}$ 
in the $(\tilde x,\sigma,\tilde z)$ coordinates.
Recall that $M^\eps_{2D}$ is strongly attracting, which means that it
inherits the stability property of $S_{a,\eps}$.
The procedure of extending slow manifolds into the flow region
is described in detail in \cite{Krupa01}.}
\end{rem}
\subsection{The manifolds $M^{\eps,r}_{2D}$ and $M^{\eps, r}_{1D}$}
\label{sec-G}
In the preceding section we constructed and analyzed
the manifold $M^\eps_{2D}$ as the continuation
of $S_{a,\eps}$ in the $(\tilde x,\sigma,\tilde z)$ coordinates.
In this section we carry out an analogous procedure, namely
we find a manifold $M^{r,\eps}_{2D}$ that extends
the unstable slow manifold $S_{r,\eps}$ backward in time.
Note that \eqref{eq-rescD1crp} has a equilibria given by the conditions 
$\sigma=0$ and $\tilde x= 1+O(\omega)$.
Let $M^{\eps,r}_{2D}$ be a center manifold of this surface of equilibria.
This manifold satisfies the condition
\begin{equation}\label{eq-condre}
\tilde x= 1+O(\sigma^2,\sov).
\end{equation}
The manifold $M^{\eps,r}_{2D}$ is foliated by  two dimensional surfaces 
$\sigma \omega=\eps$, defining invariant submanifolds of
\eqref{eq-rescD1cr}. 

We now substitute the constraint \eqref{eq-condre}
into the $(\sigma, z,\sov )$ subsystem  of
\eqref{eq-rescD1crp} and cancel a factor of $\sigma$, analogously
to the derivation of \eqref{eq-rcrcm} from 
\eqref{eq-resccr}. We obtain:
\begin{align}\label{eq-onMre2D}
\begin{split}
&\sigma '=-\frac12\sigma^2(1-\tilde z+O(\sigma,\sov))\\
&\tilde z'=\sigma\bar\mu+a(1+O(\sigma^2,\sov))+ b\tilde z
-\frac12\sigma\tilde z(1-\tilde z+O(\sigma,\sov))\\
&\sov '=\frac12\sigma\omega(1-\tilde z+O(\sigma,\sov)).
\end{split}
\end{align}
The system \eqref{eq-onMre2D} has an attracting two dimensional
invariant manifold $M^{\eps, r}_{1D}$ given by condition \eqref{eq-condre} and
\begin{align*}
\tilde z=-\frac{a}{b}+O(\sigma, \sov).
\end{align*}
The manifold $M^{\eps, r}_{1D}$ defines one dimensional submanifolds
of \eqref{eq-rescD1crp} via the constraint $\eps=\sigma\sov$.
The manifold $M^{r,\eps}_{2D}$, as defined, is not unique.
We add the  requirement that $M^{r,\eps}_{2D}$ is equal to $S_{r,\eps}$ in
the section $\tilde\Sigma^{in}$ (see Remark \ref{rem-mr2de}).
A very important feature of $M^{r,\eps}_{2D}$ is that 
$\tilde z$ is contracting forwards in time, which means that
if the manifold is continued backwards in time, the $\tilde z$ direction
is rapidly expanding.
In the $\eps\to 0$ limit there exist manifolds $M^{r}_{2D}$ and
$M^{r}_{1D}$ of \eqref{eq-resccr} that are analogous
to $M_{2D}$ and $M_{1D}$.
\subsection{A bound for the MMO behavior}
As mentioned in Section \ref{sec-micro}, MMOs cannot exist
if $\bar\mu$ is too large. In other words, the region of the
the existence of MMOs in the parameter $\mu$ must be  $O(\eps)$ in size.
Here we state and prove the relevant result. 
\begin{prop}\label{prop-bound}
There exists $L>0$ and $\eps_0>0$ such that for $\eps<\eps_0$
and $\bar\mu> L$ no MMOs can exist.
\end{prop}
{\bf Proof}
  By the analysis of Section \ref{sec-cont}, 
we know that solutions originating in $S_{1D,\eps}$,
transformed to the variables $(\tilde x,\sigma , \tilde z)$, 
satisfy $\tilde x =-1+O(\sigma^2, \omega)$ and $\tilde z=\frac{a}{b}
+O(\sigma, \omega)$. By assumption \eqref{eq-ab}, we have $\tilde z>\tilde x$
for sufficiently small $\sigma$ and $w$. By possibly decreasing the bound on
$\sigma$, we can make sure that $\tilde z>\tilde x$ is satisfied even when
we increase $\bar \mu $ (see \eqref{eq-rescD1crp}).
We transform to the $(\bar x, \bar y, \bar z)$ variables at $\sigma=\delta$,
where $\delta$ is very small and possibly depends on $\bar\mu$.
The corresponding value of $\bar y$ is $\bar y=1/\delta^2\equiv K$.
We carry out the remainder of the argument for system \eqref{eq-resc},
leaving the generalization to \eqref{eq-rescD1} to the reader.

We choose a constant $K_1$ satisfying
$K>K_1>0$ so that we can make a ``tube''
around the parabolic cylinder $\bar y=\bar x^2$, with the top and the bottom given
by $\bar y= K$ and $\bar y=K_1$ so that the solutions
can exit the tube only through the bottom. We define as sides of the tube
the parabolic cylinders
$\bar y=(1\pm\alpha)\bar x^2$, where $\alpha>0$ is a small constant, and leave the tube
unbounded in the $\bar z$ direction. 
For the moment we restrict attention
to the part of the trajectory for which $\bar z\le 0$ (and $\bar x\le 0$).
We now prove that if $K$ and $K_1$ are sufficiently
large, the trajectories
cannot exit the tube through the sides. Note that due to $\bar z\le 0$
in \eqref{eq-resc}, we obtain $\bar x'=\mp\alpha x^2$ on  
$\bar y=(1\pm \alpha)\bar x^2$ with $|\bar y'|\le|\bar x|$ in both cases.
If $K_1$ is sufficiently large, the sides of the tube are almost
vertical, and the $(\bar x, \bar y)$ component of the vector field
restricted to the sides of the tube is almost horizontal and must
be pointing towards the interior of the tube. It follows that trajectories
can only exit through the bottom, as claimed. 
By increasing $\bar \mu$ we can guarantee that
$\bar z= 0$ when the trajectory reaches the bottom of the tube.

Now note that $\bar z$ will remain non-negative. Also, 
as long as $\bar x$ remains non-positive,
$\bar y $ must decrease and
$\bar z\le -\bar \mu/b$ (since $\bar z'\le \bar \mu+b\bar z$).
However, $\bar x$ cannot reach $0$ as long as $\bar y$ is positive.
Also, since $\bar y\le K_1$ we must have $\bar x\ge -\sqrt{K_1/(1-\alpha)}$.
Hence, $\bar y'\ge -(\sqrt{K_1/(1-\alpha)}-\bar\mu/b)$. It is
now easy to estimate that the time T need for $\bar y$ to reach $0$
satisfies 
\[
T\ge \frac{K_1}{\sqrt{K_1/(1-\alpha)}-\frac{\bar\mu}{b}}.
\]
By taking $\bar \mu $ sufficiently large we can guarantee that 
by the time $\bar y$ reaches $0$, $\bar z$ must
become at least $-bK_1/2$. By increasing $\bar\mu$ we can also make 
$K_1$ arbitrarily large. Since $\bar z$ will not decrease, 
it follows that $\bar x$ must become at least equal to
$-bK_1/2$ before $\bar y$ starts increasing. It is well known that
for any initial condition $(\bar x, \bar y, \bar z)$ with
$\bar y\le 0$, $\bar z\ge 0$, and $\bar x$ sufficiently large
the solution must enter a relaxation loop \cite{Krupa01}.
Hence, no STO (which requires $\bar y$ to increase and $\bar x$ to decrease)
can exist.\qed

\subsection{Canard solutions and the boundaries of MMO behavior}
\label{sec-pth} 
In this section we discuss the role of canards in transitions between
different MMO regimes. We explain in more detail what we mean by MMOs
and canards and discuss how transitions between different MMO regimes
can occur. We identify one type of the possible transitions as occurring
through a canard and present numerical evidence showing that the final 
transition to spiking must be of this
type. Finally we discuss the $\eps\to 0$ limit of canards in the
$(\bar x,\bar y, \bar z)$ variables and show that they
correspond to {\em singular canards} that are special solutions of 
\eqref{eq-resc}. An explicit example of such a special solution is
the algebraic solution $\gamma_{poly}$ found is Section \ref{sec-spec}.
We conjecture that the $\bar\mu$ value corresponding to $\gamma_{poly}$
is the $\eps=0$ limit of the $\bar\mu$ values where the transition between MMOs
and spiking occurs.

We begin by giving a more precise meaning to MMOs and canards.
Let $V$ be the interior of the cube 
$[-\delta, \delta]^3$, where $\delta\ge\rho^2$ 
is some fixed small number. Consider a trajectory $(x(t),y(t),z(t))$
originating near $S_{a,\eps}$ outside of V
(with the initial value of $y\ge\delta$). 
We say that $(x(t),y(t),z(t))$ has 
{\em rotation type $k$} if $x(t)$ has $k$ non-degenerate maxima and minima
during 
the passage of the trajectory through $V$.  
A trajectory originating in $S_{a,\eps}$ is a {\em $k$th secondary canard} if
$x(t)$ has $k$ non-degenerate maxima and minima during the passage through 
$V$ and subsequently continues into $S_{r,\eps}$. 
A {\em primary canard} is a trajectory that passes through $V$ with $x(t)$
monotonically increasing ($x'(t)>0$) and subsequently continues 
into $S_{r,\eps}$. 

Before we discuss MMOs and their transitions, we 
point out a difference between \eqref{canon3dlinear4}
and the folded node problem \cite{Benoit90, Szmolyan01, Wechselberger05}.
Recall the one dimensional invariant manifold $S_{1D,\eps}$ 
contained in $S_{a,\eps}$. It was 
shown in Section \ref{sec-S1D} that all trajectories on $S_{a,\eps}$ 
are attracted to $S_{1D,\eps}$.
We define an $S_{1D}$ {\em canard} as a trajectory originating in 
$S_{1D,\eps}$ and continuing to $S_{r,\eps}$. 
Note that an $S_{1D}$ canard corresponds to an intersection of a 1D
invariant manifold $S_{1D,\eps}$ and 2D invariant manifold $S_{r,\eps}$.
If such an intersection is transverse with respect to the variation 
of $\bar\mu$ (i.e. in the 4D extended phase space including
$\bar\mu$), then $S_{1D,\eps}$ crosses from one
to the other side of $S_{r,\eps}$. If we knew that all the intersections
were transverse, we could conclude that $S_{1D}$ canards exist only
for isolated values of $\bar\mu$. 
Suppose that $\gamma(t)$ is a canard solution that 
is not an $S_{1D}$ canard (the initial condition is not contained in
$S_{1D,\eps}$). Since $\gamma(t)$, during
its passage along $S_{a,\eps}$, becomes exponentially close 
to $S_{1D,\eps}$, it follows that the (minimal) separation between 
$S_{1D,\eps}$ and $S_{r,\eps}$ must be exponentially small. 
Consequently, if we knew that the distance between $S_{1D,\eps}$
and $S_{r,\eps}$ as a function of $\bar\mu$
were not identically equal to 0 and varied at least algebraically in
$\bar \mu$, we could conclude that
canard trajectories can only exist on exponentially small 
($O(e^{-c/\eps}$), $c=\mbox{ const}>0$) intervals of the parameter $\bar\mu$.
By contrast, in the folded node problem, canards exist robustly (for all the
values of the control parameter) \cite{Benoit90, Szmolyan01, Wechselberger05}.

We now consider transitions between different types of MMOs.
We will restrict attention to MMO periodic orbits of type $1^k$,
i.e. $k$ small loops followed by a spike. More precisely, we consider 
MMO periodic orbits with rotation type $k$ (see the definition above).
If we change $\bar \mu$ then the rotation type can change in two ways:
\begin{description}
\item[(a)]
Either a large loop can decrease in size and become fully contained in 
$V$, or a loop contained in $V$ can grow and be no longer contained in $V$.
\item[(b)]
For some intermediate value of $\bar\mu$, $x(t)$ has 
a degenerate critical point corresponding to the collapse of a small loop.
\end{description} 
The following result shows that if a certain non-degeneracy condition 
holds, then all the transitions are of type (a).  Moreover, 
they correspond to transitions through a canard.
\begin{prop}\label{prop-tran}
Suppose there exist two MMO periodic orbits with rotation type $k$
and $k-1$, respectively, for two different values of $\bar\mu$, denoted by
$\bar\mu_1$ and $\bar\mu_2$, respectively. Recall the special
solution $\gamma_{\bar\mu}=(\bar x_{\bar\mu}(t),\bar y_{\bar\mu}(t),
\bar x_{\bar\mu}(t))$ and assume that for any $\bar\mu\in
[\bar\mu_1,\bar\mu_2]$ $\bar x_{\bar\mu}(t)$ has only non-degenerate critical
points. Then, provided that $\eps$ is sufficiently small,
the transition between the two MMOs is through a canard.
\end{prop}
{\bf Proof}
Let $\gamma_1$ and $\gamma_2$ denote the MMO with rotation type $k$
and $k-1$, respectively. Suppose we vary $\bar\mu$ from $\bar\mu_1$
to $\bar\mu_2$. Given a trajectory $(x(t), y(t), z(t))$ of \eqref{canon3d}, 
note that a critical point of $\bar x(t)$ corresponds
to a crossing with the nullsurface $y= x^2+D_1x^3$.
At such a critical point we have $x''=-y'=-\eps^2(x-z)$, so that
a degenerate critical point can only occur if $x=z$. If the
solution crosses the nullsurface $y= x^2+D_1x^3$
outside of $V$ and close to $S_{r,\eps}$, then the solution
must have followed $S_{r,\eps}$ for a significant amount of time
so that $\mu+ax+bz=O(\eps)$. Thus $ax+bz=O(\eps)$, 
by Proposition \ref{prop-bound}. 
Hence, for $\eps$ sufficiently small, the critical point
must be non-degenerate. 
Now we consider a crossing of the x-nullsurface  inside of V.
Any such crossing point must also be non-degenerate,
due to the assumption on $\gamma_{\bar\mu}$ (by compactness
and smooth dependence on initial conditions, the intersections
of $M_{1D}^\eps$ with $\bar y=\bar x^2+\eps D_1 \bar x^3$
must also be non-degenerate). 
For $(\bar x, \bar y)$ outside of a bounded region
we switch to the $(\tilde x,\sigma,\tilde z,\omega)$
coordinates and observe that, in order to reach the nullsurface, 
the trajectory would have to follow $M_{2D}^{r,\eps}$, 
and therefore be attracted to $M_{1D}^{r,\eps}$, so that
$\tilde z\approx -a/b$ would have to hold (see Section \ref{sec-G}), 
making $x=z$ impossible since $a<-b$.
We have proved that all the crossing points with the $y= x^2+D_1x^3$
nullsurface must be non-degenerate. It follows that the only way a transition
from $\gamma_1$ to $\gamma_2$ can occur is if a point of intersection
of the trajectory with the nullsurface $y= x^2+D_1x^3$ traces the 
entire nullsurface until one of the small loops develops to a spike. \qed

\begin{figure}[!htb] 
\begin{center}
\includegraphics[width=.45\textwidth]{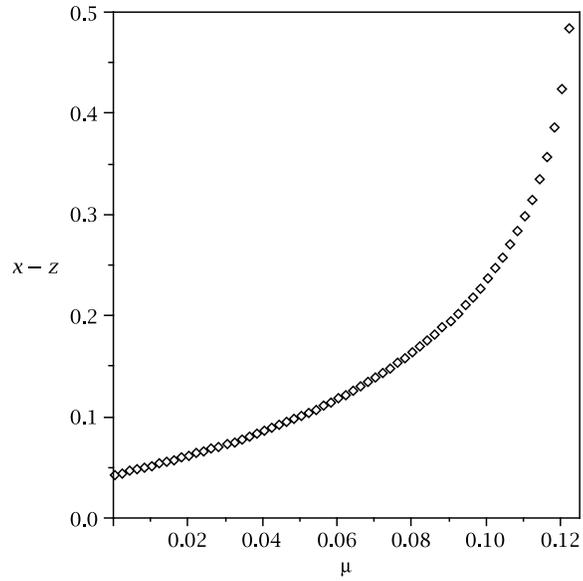}
\caption{Plot of  of $\bar x -\bar z$ vs. $\bar\mu$ at the first intersection point of $\bar\gamma_{\bar\mu}$
with the surface $\bar y=\bar x^2$ in \eqref{eq-resc}.}\label{tran1fig}
\end{center} 
\end{figure}

Although we believe that all the intersections of $\gamma_{\bar\mu}$
with the nullcline $\bar y= \bar x^2$ (for system \eqref{eq-resc}
are transverse, we have only
partial results in this direction, and we feel that they are not suitable for
this article and will be presented in future work. Here we state the
following conjecture.
\begin{conj}\label{conj-tran}
The least $t$ for which $\bar x'_{\bar\mu}= 0$ (for system \eqref{eq-resc})
corresponds to a transverse intersection.
\end{conj}
As support for this conjecture we have plotted the value of
$\bar x -\bar z$ at the first intersection point of $\bar\gamma_{\bar\mu}$
with the surface $\bar y=\bar x^2$ in Figure \ref{tran1fig}.
\begin{rem}\label{rem-mmocan}
Suppose that conjecture \ref{conj-tran} is true. Then modifying slightly the
proof of Proposition \ref{prop-tran} we could prove that the transition
from MMO to spiking is through a primary canard.
\end{rem}
The following result concerns the $\eps\to 0$ limit of canards in the
context of \eqref{eq-rescD1}. We say that the special
solution $\gamma_{\bar\mu}$
\eqref{eq-resc} is a {\em singular canard} if it continues to $M^r_{2D}$.
\begin{thm}\label{thm-main}
Suppose $\bar\mu(\eps)$ is a continuous function 
such that for each $(\bar\mu,\eps)$ satisfying
$\bar\mu=\bar\mu(\eps)$ system \eqref{eq-rescD1} 
has a primary canard.  
Then, in the $\eps\to 0$ limit, 
system \eqref{eq-resc}  has a singular canard. 
\end{thm}
{\bf Proof}
Recall that in Section \ref{sec-cont} we showed that the continuation 
of $S_{1D,\eps}$ to the fold region must be $O(\eps)$ close to 
the special solution $\gamma_{\bar\mu}$. As shown in Section \ref{sec-G},
equation \eqref{eq-rescD1c}
has a 2-d repelling invariant manifold $M_{2D}^{r,\eps}$, 
which is a continuation of $S_{r,\eps}$ to the fold region (backward in time). 
In the limit $\eps\to 0$ we obtain a two dimensional invariant  manifold 
of \eqref{eq-rescc}, which we denote by $M_{2D}^r$.  
A canard solution, as it passes through the fold region, 
must become exponentially close to $M_{1D}^\eps$, 
and then connect to $M_{2D}^{r,\eps}$, before going
on to $S_{r,\eps}$. In the $\eps\to 0$ limit we get a connection 
from $\gamma_{\bar\mu}$ to $M_{2D}^r$. The $\bar\mu$ value corresponding 
to an $\eps\to 0$ limit of a transition from MMOs to spiking has to be one 
for which such a connection exists. Note that the validity
of this argument relies on the convergence of the continuation 
of $S_{1D,\eps}$ to the special solution 
$\gamma_{\bar\mu}$ in the variables $(\tilde x, \sigma, \tilde z)$.
The relevant convergence result was established in Section \ref{sec-cont}.\qed
\begin{rem}\label{rem-spec}
The defining property of a singular canard, namely $\gamma_{\bar\mu}$ 
continuing to $M_{2D}^r$, is non-generic and is likely to occur for
isolated values of $\bar\mu$ (provided that a suitable transversality
condition holds). We know one example, 
with $\bar\mu=\bar\mu_0=-\frac{a(a+b)}{2b^2}$, corresponding to 
the polynomial vector solution $\gamma_{{\rm poly}}=\gamma_{\bar\mu_0}$.
We conjecture that $\bar\mu_0$ is the $\eps\to 0$
limit of the locus of the transition from MMO to spiking. This conjecture
will be presented in more detail at the beginning of the next section.
\end{rem}
\subsection{Transition from MMO to spiking -- numerical results} 
\label{sec-percan}

We propose that $\gamma_{poly}$ is the unique singular canard
and that it persists as a regular $S_{1D}$ canard, that is
a connection between $M_{1D}^\eps$ and $M_{2D}^{r,\eps}$.
Consequently, we conjecture the following: 
\begin{conj}\label{conj-main}
There exists a function $\bar\mu (\eps)$ such that for $\mu=\mu(\eps)$
system \eqref{canon3d} has a primary $S_{1D}$ canard. The transition 
from MMOs to spiking occurs within an exponentially small ($O(e^{-c/\eps})$, $c=\mbox{ const
 }>0$) interval of $\bar\mu$ about $\bar\mu(\eps)$. Moreover,
\begin{align*}
\bar\mu(0)=\bar\mu_0=-\frac{a(a+b)}{2b^2}.
\end{align*}
\end{conj}

\begin{figure}[!htb] 
\begin{center}
\includegraphics[width=.8\textwidth]{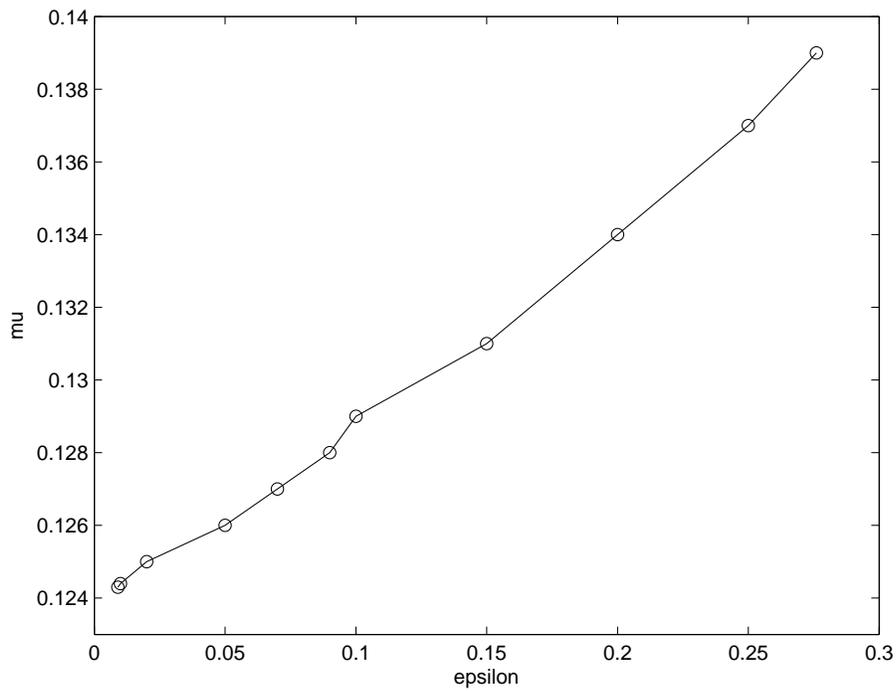}
\caption{Values of $\bar \mu$ corresponding to the transition between
MMOs and pure spiking in system \eqref{canon3d} for decreasing values of epsilon.}\label{muepsfig}
\end{center} 
\end{figure}

We investigate this conjecture numerically by first considering the transition threshold between MMOs and pure spiking  in system \eqref{eq-rescD1} by computing solution trajectories for decreasing values of $\eps>0$ and $\bar\mu$ near $\bar\mu_0$  The simulations show that for decreasing values of $\eps>0$, the transition between MMOs and pure spiking occurs for decreasing values of $\bar\mu$ and suggest that $\bar\mu \to \bar\mu(0)\approx 0.1238928363$ as $\eps \to 0$ as seen in Figure \ref{muepsfig}. 

\begin{figure}[!htb] 
\begin{center}
\includegraphics[angle=270, width=\textwidth]{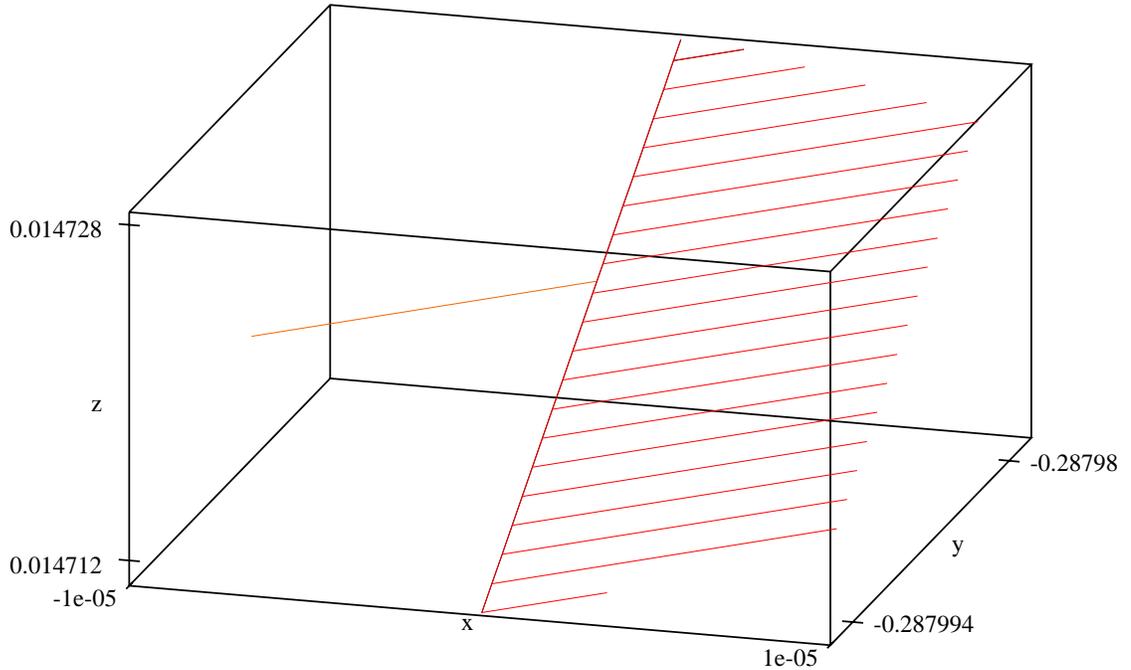}
\caption{Approximation of the primary canard formed by the connection of approximations of $M_{1D}^\eps$ and $M_{2D}^{r,\eps}$ for $\bar\mu=.1381943$ in system \eqref{eq-rescD1} with the default value of $\eps=.2764436178$.  The family of curves approximating the manifold $M_{2D}^{r,\eps}$ is obtained by taking a sequence of initial conditions with $\bar x=5$, $\bar y= 10$ and $\bar z$ near $1.917386$ and evolving backwards in time until reaching the $\bar x=0$ plane. The single curve approximating $M_{1D}^\eps$ is generated by evolving in forwards time until reaching the $\bar x=0$ plane with initial conditions $\bar x=-2$, $\bar y=3.726472$, $\bar z=-.9058866$ on the special solution $\gamma_{{\rm poly}}(t)$ with $\bar y$ and $\bar z$ as prescribed in Section \ref{sec-spec}.}\label{connect3d}
\end{center} 
\end{figure}

\begin{figure}[!htb] 
\begin{center}
\includegraphics[angle=270, width=.5\textwidth]{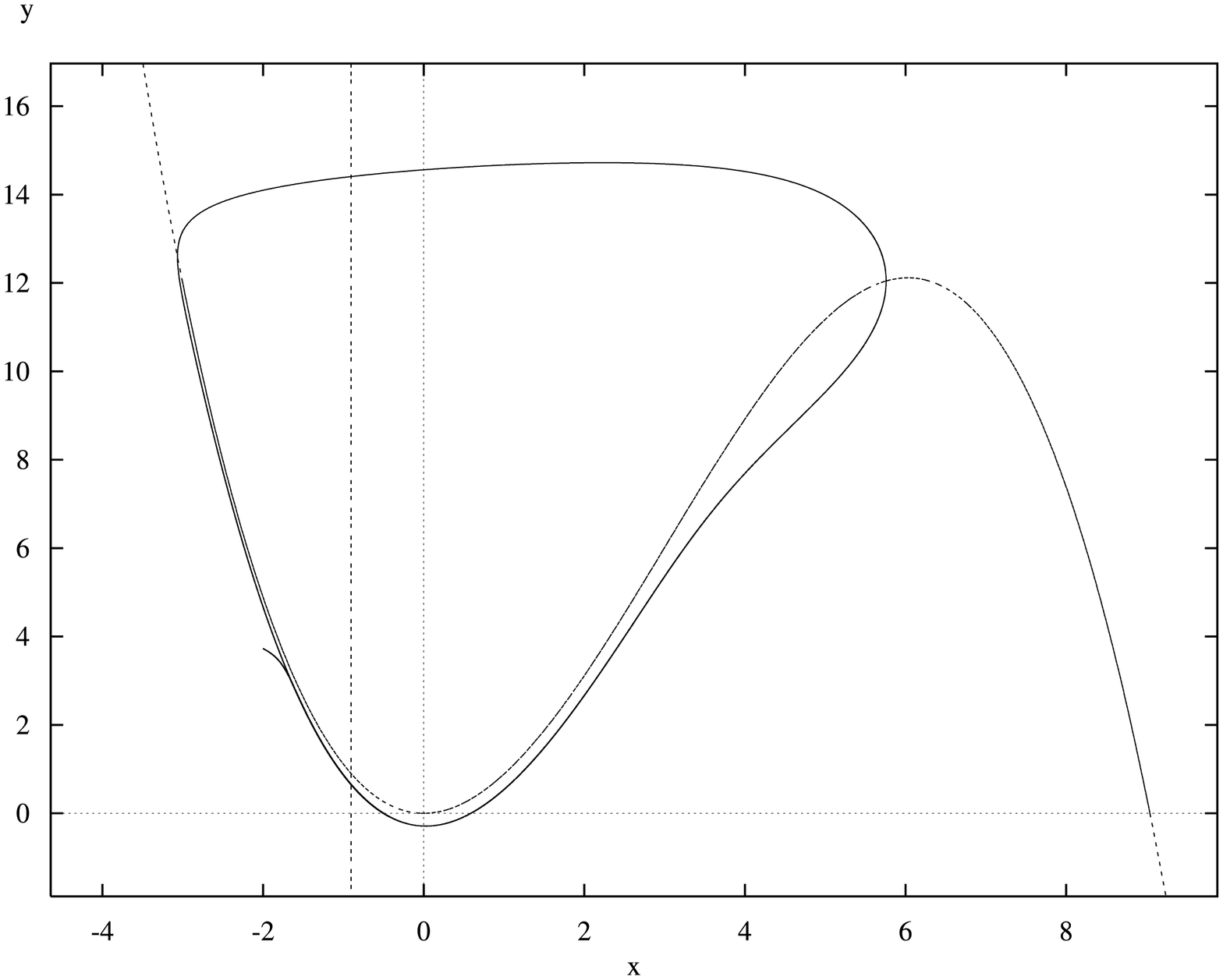}
\includegraphics[angle=270, width=.5\textwidth]{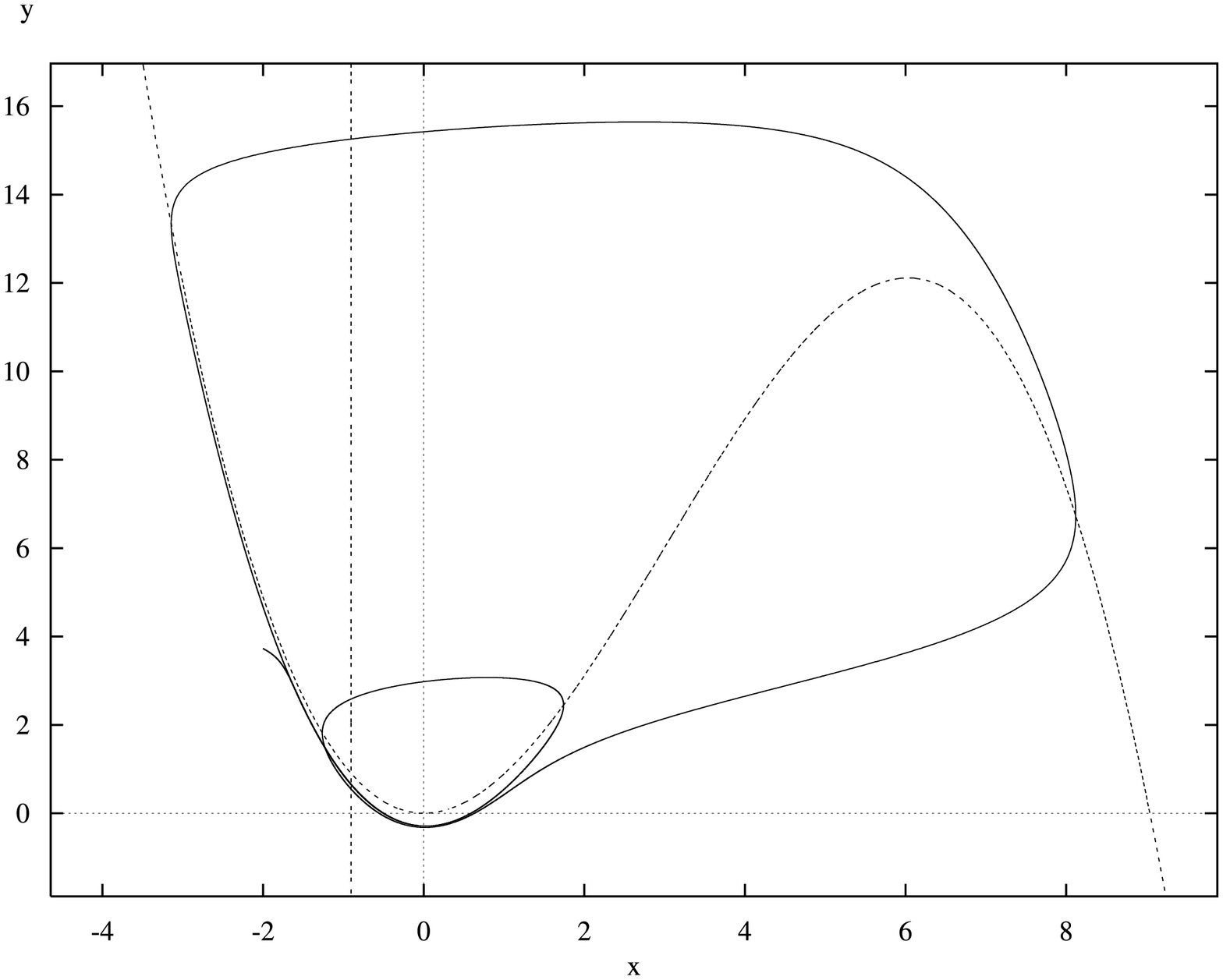}
\includegraphics[angle=270, width=.5\textwidth]{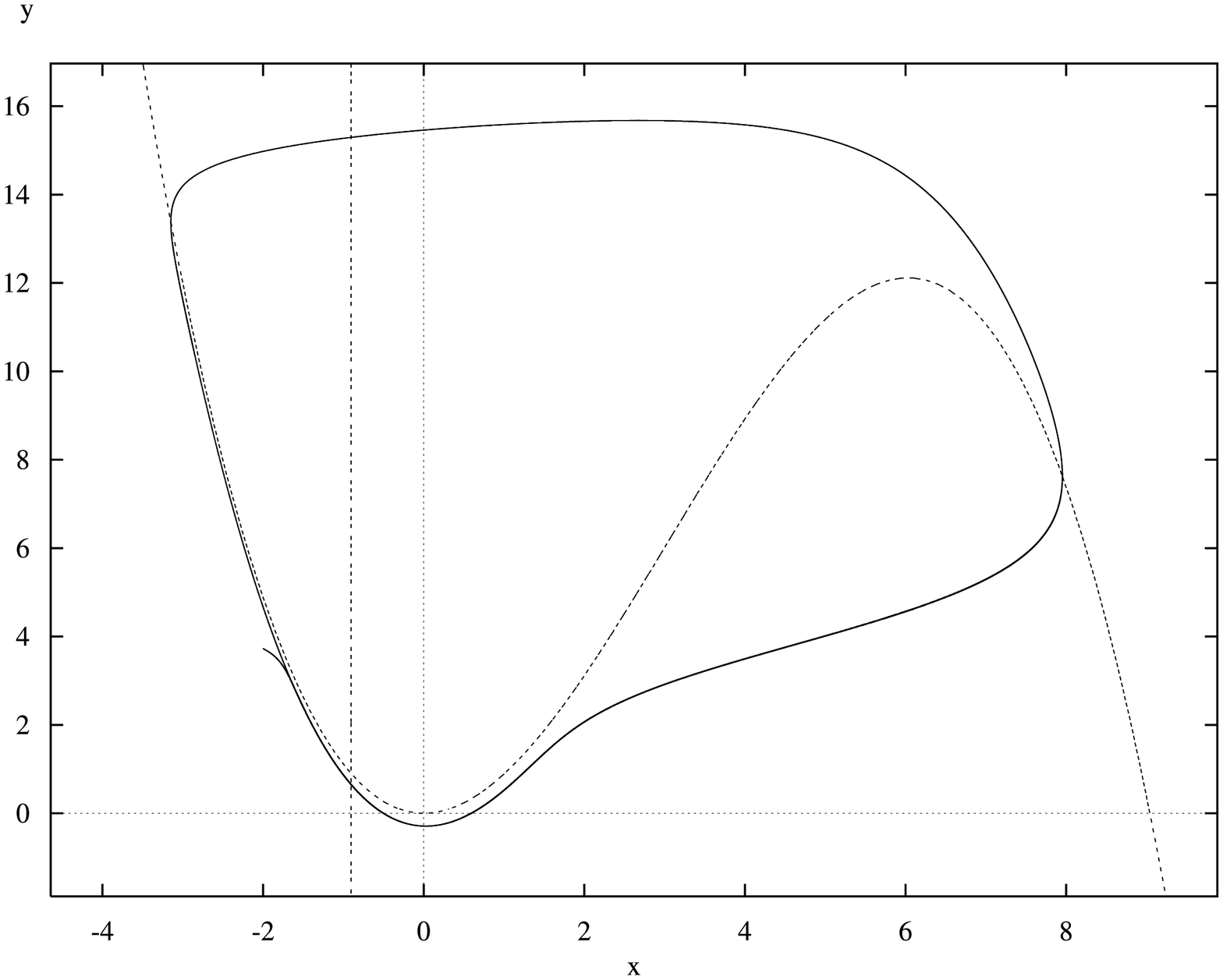}
\caption{Top: Approximate primary canard solution of system \eqref{eq-rescD1} with $\bar\mu=0.1381943$ near the transition between mixed-mode oscillations and spiking. Middle: Approximate primary canard solution of system \eqref{eq-rescD1} with $\bar\mu=0.1381943$ near the transition between mixed-mode oscillations and spiking. Bottom: Periodically spiking solutions exhibited by system \eqref{eq-rescD1} with $\bar\mu=0.139$. }\label{3traj}
\end{center} 
\end{figure}

In addition, we have found an approximate primary canard numerically.  As noted above, such a solution must become exponentially close to the attracting manifold $M_{1D}^\eps$ as it passes through the fold region and then connect to the repelling manifold $M_{2D}^{r,\eps}$ before going on to $S_{r,\eps}$.  We exhibit this numerically using the Poincare section option of XPPAUT (Ermentrout02) for $\eps=\eps_d$.  That is, we choose an initial condition with $\bar x=-2$ on the singular primary canard $\gamma_{{\rm poly}}$ and integrate until reaching the Poincare section $\bar x=0$.  See Figure \ref{connect3d}. This solution is expected to become exponentially close to $M_{1D}^\eps$.  Next, we approximate $M_{2D}^{r,\eps}$ by choosing a range of initial conditions varying in $\bar z$ with $\bar x=5$, $\bar y= 10$, and $\bar z$ near $1.917386$ and integrating backwards in time up to the Poincare section at $\bar x=0$ to form a family of curves which form an approximation of $M_{2D}^{r,\eps}$.  As shown in Figure \ref{connect3d}, an approximate connection between these manifolds is obtained for $\bar\mu=.1381943$.  The approximate primary canard solution that corresponds to the intersection of these two manifolds and shows the characteristic long term (O(1)) following of a repelling manifold is also displayed in Figure \ref{3traj}. 

Additionally, numerical simulations with varying levels of  $\bar\mu$ show a transition between MMOs and spiking near the value of $\bar\mu$ for the approximate primary canard as seen by the MMO solution trajectory for $\bar\mu=.138$ and the spiking solution for $\bar\mu=.139$ in Figure \ref{3traj}.
Thus, our numerical simulations support the validity of the persistence of the primary canard as stated in Conjecture \ref{conj-main}.

There are many other interesting issues to analyze in this model including the transition from small amplitude dynamics to MMOs and the rotation properties of MMOs.  For example, one could try to obtain an approximation for the number of subthreshold oscillations per spike for a given value of $\bar\mu$. We discuss such issues in the Discussion, but a detailed analysis of these issues is beyond the scope of this paper.

\section{Discussion}\label{sec-discussion}

In this work, we have analyzed a model introduced by Acker et al.\cite{Acker03} 
to study synchronization properties of stellate cells in layer II of the entorhinal cortex (EC).    We were interested in studying subthreshold oscillations in pyramidal cells of layer V of the EC.  Since there was no available model of EC layer V pyramidal cells at the time, and since the two cell types exhibited very similar electrophysiological characteristics, including the presence of persistent sodium and slow potassium currents that are thought to underlie the dynamics of STOs, we chose to employ the Acker model to model STOs in layer V pyramids. 
We note that layer V pyramidal cells exhibit many other interesting dynamical properties.
 They participate in high frequency ($\sim 200$ Hz) "ripple" activity associated with hippocampal sharp waves \cite{Chrobak94} and are a site of initiation and propagation of epileptiform activity\cite{Jones90}.  Also, they exhibit graded persistent activity, which is associated with memory processes, in which neurons respond to consecutive stimuli with graded changes in firing frequency that remain stable after each stimulus \cite{Egorov02}.  In 2006, Fransen et al. introduced a quite detailed multi-compartmental model of layer V pyramidal cells that exhibit the graded persistent activity.  In future work, we plan to adapt our analysis of the simpler model considered here to this more realistic model, and analyze the dynamics of the graded persistent activity. 
 
We began our work with a numerical investigation of the dynamical properties and bifurcation structure of the full model and displayed the role of the persistent sodium and slow potassium currents in the generation of subthreshold oscillations.  In order to analyze the 6-d system, we performed a reduction of dimension analysis to reduce the system to a 3-d system. We note that a more thorough reduction of dimensions analysis of a similar model of a layer II stellate cells with an h-current and persistent sodium current was performed by Rotstein et al. \cite{Rotstein06}, and a computational tool for reducing systems possessing multiple scales has also been developed \cite{kn:clekop1}.  
Rotstein et al. have also analyzed MMOs present in that model \cite{rotsteinwechselberger07}. 

Just as experimental neuroscientists often modulate or block various currents in order to learn about their role by observing the effects of their absence, we investigated the effect of modulating various currents in our model.  Upon removing or reducing to a very low level the persistent sodium current and replacing it with a tonic excitatory input, we witnessed very interesting dynamical phenomena including MMOs and seemingly chaotic firing patterns.  (Future work includes the investigation of whether our analysis and the dynamical mechanism involved carry over to the case where a small persistent sodium current is retained.)  This experiment is unable to be performed presently in the lab since there are no known drugs which are able to selectively block persistent sodium current without blocking the standard sodium current, and we felt it would be interesting to investigate this regime, especially due to the rich dynamical behavior present.  The main focus of this work is to explain the dynamical mechanisms responsible for the mixed-mode oscillations.

We have discovered and analyzed a novel 
dynamical mechanism for the generation of MMOs,
involving the canard phenomenon
as well as the dynamics near a saddle focus equilibrium, as in the 
Shilnikov scenario \cite{Koper92} or the subcritical Hopf scenario
\cite{Guckenheimer1, Guckenheimer2}.  
The techniques we employed to analyze the MMOs involve standard techniques from dynamical systems theory such as geometric singular perturbation theory \cite{Fenichel79,Jones95,Krupa01} and center manifold theory \cite{Carr81}. For a summary of the analysis, see the last paragraph of section \ref{sec-intro}. 
 
A main accomplishment of this work is an analytical expression for the parameter at which the transition between MMOs and spiking occurs in the singular limit and a characterization of the corresponding singular primary canard solution responsible for the transition.  Moreover, we construct a well developed geometric picture of the dynamics for this new canard-based MMO mechanism since we utilize a great deal of the underlying geometric structure available including slow manifolds, canards, and the saddle-focus equilibrium with its weak 
unstable manifold.  In addition, the use of the folded-node point from the two time scale model to derive our novel three time scale model is new. 
Also, our analysis, which resulted from the investigation of a neuronal model, is applicable to studying MMOs in other systems from various contexts that have the same underlying dynamical structure, i.e. systems that can be transformed into the form of \eqref{canon3d}.

Prominent among scenarios for MMOs
not involving canards are the Shilnikov scenario
\cite{Shilnikov65,Koper95,Rajesh00} 
and the subcritical Hopf scenario \cite{Guckenheimer1, Guckenheimer2}.  
Both involve a saddle focus equilibrium with a two dimensional
unstable manifold and a return mechanism
of 'soft reinjection', which basically means that the return mechanism is given by a typical orientation preserving diffeomorphism.
Our model system \eqref{canon3d} has a saddle focus equilibrium with
a weak unstable manifold, but in addition it has a multiple time scale
structure, folded critical manifolds, and a singularity similar to a folded
node. One of the consequences is that the return mechanism is not soft, but
has a definite structure determined by the slow manifolds. Moreover,
as we have demonstrated in the paper, the multiple time scale structure
influences subthreshold oscillations by means of canard solutions.

It follows from our analysis that all the orbits leaving the vicinity
of the saddle focus equilibrium will be attracted to the super-stable
slow manifold, which returns to the vicinity of the equilibrium.
Thus the system can be viewed as close to a homoclinic bifurcation. 
In this paper we have concentrated on the aspects of the
problem related to the canard mechanism,
 but we feel that the methods employed to study homoclinic  
bifurcations, and especially the passage near a saddle point,
 may be applicable and useful
for a more detailed investigation of our problem.
We intend to explore the homoclinic-like features of the problem 
in future work.

In future work, we also intend to analyze the rotational properties of the solutions of \eqref{canon3d}.
For a solution originating in $S_{a,\eps}$ or nearby, we define
its {\em rotation number} to be the number of small loops it makes while
passing through the {\em fold region}, which we define as a sufficiently
large ball centered at the origin 
in the coordinates $(\bar x,\bar y, \bar z)$. Note that due to the strong
attraction of $S_{1D,\eps}$ , all the solutions 
entering the fold region  
will have the same rotation number,
unless the system is very close to a canard. Here we mean a primary canard
or a secondary canard, as secondary canards correspond to boundary curves, in 
the $(\mu,\eps)$ parameter space, between parameter
regions for which solutions have different rotation numbers.
We conjecture that for a fixed $\eps$, there is a sequence of $\bar\mu$ values,
corresponding to $k$ secondary canards with increasing $k$
and to transitions corresponding to an increment by $1$  in the
rotational number, ending with a region corresponding to
the maximal number of rotations beyond which only small amplitude dynamics
exists.  Also, we have not analyzed the transition from small amplitude dynamics
to MMOs. We know from simulations that small amplitude chaos is present,
and we conjecture that it is related to the loss of normal hyperbolicity
of the unstable manifold of the saddle-focus equilibrium,
which must occur sufficiently far away from the equilibrium. Also, we believe that it is related to the Shilnikov scenario involving incomplete homoclinicity.  

Another interesting question is how to estimate the rotation number for
given values of $\bar\mu$ and $\eps$. Unfortunately, we do not foresee that
an analytic formula can be obtained, in contrast to related cases in which there are two timescales \cite{Szmolyan01}. 
A conjecture very strongly supported by the numerics is that 
the rotations happen almost exclusively during the slow drift of the solutions
along the unstable manifold of of the saddle focus equilibrium.
We also believe that the flow along this unstable manifold is close to 
its linear approximation. Once the solutions  come to a region where 
$\bar x$ and/or $\bar y$ are sufficiently large, the unstable manifold 
loses normal hyperbolicity and the multiple time scale nature of the system
takes over, thus ending the rotational phase of the solution.
However, in order to produce an approximate formula for the number of
rotations, one would need to know the distance from the continuation of
$S_{1D,\eps}$ to the stable manifold of the saddle focus equilibrium, measured in
a specified section such as the hyperplane
$\bar x=0$, as well as the distance between the
equilibrium and the region where normal hyperbolicity is lost. One would 
also need the information on the eigenvalues and the eigenvectors of the
saddle focus equilibrium. 
No analytic version of the formula can be expected since it seems quite unlikely 
to obtain an analytic formula for either the continuation of $S_{1D,\eps}$
all the way to the section $\bar x=0$ or the stable manifold of the saddle
focus equilibrium, but a numerical formula seems achievable.

Note that even in the case of \cite{Szmolyan01}, 
only a formula asymptotic 
in $\eps$ is available, and in fact the formula of Szmolyan and Wechselberger 
is the $\eps=0$ approximation. In our case, we can consider an asymptotic
formula based on the linear approximation of the
flow along the unstable manifold,
and in the $\eps=0$ approximation this would mean that $\gamma_{\bar\mu}$
would be used instead of
the continuation of $S_{1D,\eps}$. However, again no analytic formula can be expected
as no analytic expression for $\gamma_{\bar\mu}$ is available, except for the special case where $\bar\mu=\bar\mu_0$.
A possibility would be to obtain an asymptotic expression using
center manifold theory, along the ideas of Section \ref{sec-spec},
and we plan to pursue this idea in future work.

\clearpage
\section{Appendix} \label{section appendix}
\subsection{Model equations}
The full model for the layer V pyramidal cell is given by
\begin{eqnarray}
C\frac{dv}{dt}&=&-g_{Na}m^3h(v-E_{Na})-g_K n^4(v-E_K)- g_L(v-E_L)\nn\\
 && -g_{Nap}p(v-E_{Na})-g_{Ks}w(v-E_K)+ I_{app}  \nn\\
\frac{dm}{dt}&=& \frac{m_{\infty}(v)-m}{\tau_m(v)}\nn\\
\frac{dn}{dt}&=& \frac{n_{\infty}(v)-n}{\tau_n(v)}\\
\frac{dh}{dt}&=& \frac{h_{\infty}(v)-h}{\tau_h(v)}\nn\\
\frac{dp}{dt}&=& \frac{p_{\infty}(v)-p}{\tau_p(v)}\nn\\
\frac{dw}{dt}&=& \frac{w_{\infty}(v)-w}{\tau_w}\nn\
\end{eqnarray}

The default parameters for the system are: $C=1.5,g_{Nap}=.21$, $g_{Na}=52$, $g_K=11$, $g_L=.1$, 
$g_{Ks}=2.0$, $\tau_w=90$, $E_{Na}=55$, $E_K=-90$, and $E_L=-54$.
The graphs of the voltage dependent activation and inactivation functions as well as the voltage dependent time constants are given in Figures \ref{figactivCurves} and \ref{figRatecurves}, and 
their equations are the following: \\
$m_{\infty}(v)=\alpha_m(v)/(\alpha_m(v)+\beta_m(v))$ where \\
$\alpha_m(v)=-.1(v+23)/(\exp(-.1(v+23))-1)$ and
$\beta_m(v)=4\exp(-(v+48)/18)$.\\
$n_{\infty}(v)=\alpha_n(v)/(\alpha_n(v)+\beta_n(v))$ where\\ 
$\alpha_n(v)=-.01(v+27)/(\exp(-.1(v+27))-1)$ and
$\beta_n(v)=.125\exp(-(v+37)/80)$.\\
$h_{\infty}(v)=\alpha_h(v)/(\alpha_h(v)+\beta_h(v))$ where\\ 
$\alpha_h(v)=.07\exp(-(v+37)/20)$ and
$\beta_h(v)=1/(\exp(-.1(v+7))+1)$.\\
$p_{\infty}(v)=1/(1+\exp(-(v+38)/6.5))$.\\
$w_{\infty}(v)=1/(1+\exp(-(v+35)/6.5))$.\\
The equations for the voltage dependent time constant for $n$ is:\\
$\tau_n(v)=1/(\alpha_n(v)+\beta_n(v))$.

\clearpage
\section{Figures}\label{section figures}

\clearpage

\clearpage

\clearpage

\clearpage

\clearpage

\begin{acknowledgments}
The authors are deeply indebted to Nancy Kopell for her numerous
contributions to the paper.
The authors would also like to thank Nicola Popovi\'c, Martin Wechselberger,
and the reviewers for fruitful discussions. 
The research of MK was supported in part by NSF grant DMS-0406608.
The research of JJ was supported in part by a YSU Research Professorship.
\end{acknowledgments}

\bibliographystyle{unsrt}

\begin{thebibliography}{10}

\bibitem{Zhabotinsky64}
A.M Zhabotinsky.
\newblock Periodic kinetics of oxidation of malonic acid in solution.
\newblock {\em Biofizika}, 9:306--311, 1964.

\bibitem{Koper95}
M.T.M. Koper.
\newblock Bifurcations of mixed-mode oscillations in a three-variable
  autonomous van der pol-duffing model with a cross-shaped phase diagram.
\newblock {\em Physica D}, 80:72--94, 1995.

\bibitem{kn:milszm2}
A.~Milik, P.~Szmolyan, H.~Loffelmann, and E.~Groller.
\newblock Geometry of mixed-mode oscillations in the 3d-autocatalator.
\newblock {\em Int. J. Bif. Chaos}, 8:505--519, 1998.

\bibitem{kn:epssho1}
I.~R. Epstein and K.~Showalter.
\newblock Nonlinear chemical dynamics: oscillations, patterns and chaos.
\newblock {\em J. Phys. Chem.}, 100:13132--13147, 1996.

\bibitem{Rotstein06b}
H.~G. Rotstein and R.~Kuske.
\newblock Localized and asynchronous patterns via canards in coupled calcium
  oscillators.
\newblock {\em Physica D}, 215:46--61, 2006.

\bibitem{Chay95}
T.R. Chay, Y.S. Fan, and Y.S. Lee.
\newblock Bursting, spiking, chaos, fractals, and universality in biological
  rhythms.
\newblock {\em Int. J. Bifurc. Chaos}, 5:595, 1995.

\bibitem{DicksonOscAct00}
C.T. Dickson, Magistretti J., Shalinsky M., Hamam M., and Alonso A.
\newblock Oscillatory activity in entorhinal neurons and circuits.
\newblock {\em Annals New York Academy of Sciences}, 911:127 -- 150, 2000.

\bibitem{Yoshida07}
M.~Yoshida and Alonso A.
\newblock Cell-type specific modulation of intrinsic firing properties and
  subthreshold membrane oscillations by the m(kv7)-current in neurons of the
  entorhinal cortex.
\newblock {\em J. Neurophysiol.}, accepted.

\bibitem{medvedev04}
G.~Medvedev and J.~Cisternas.
\newblock Multimodal regimes in a compartmental model of the dopamine neuron.
\newblock {\em Physica D}, 194:333--356, 2004.

\bibitem{Drover04}
J.~Drover, J.~Rubin, J~. Su, and B.~Ermentrout.
\newblock Analysis of a canard mechanism by with excitatory synaptic coupling
  can synchronize neurons at low firing frequencies.
\newblock {\em SIAM J. Appl. Math}, 65:69 -- 92, 2004.

\bibitem{wechselbergerrubin07}
J.~Rubin and M.~Wechselberger.
\newblock The selection of mixed-mode oscillations in a Hodgkin-Huxley model
  with multiple timescales.
\newblock {\em Chaos}, 2007.

\bibitem{Rotstein06}
H.~G. Rotstein, J.A. Oppermann, T.and~White, and N.~Kopell.
\newblock The dynamic structure underlying subthreshold activity and the onset
  of spikes in a model of medial entorhinal cortex stellate cells.
\newblock {\em J. Comput Neurosci.}, 21(3):271--92, 2006.

\bibitem{rotsteinwechselberger07}
H.G. Rotstein, M.~Wechselberger, and N.~Kopell.
\newblock Canard induced mixed-mode oscillations in a medial entorhinal cortex
  layer ii stellate cell model.
\newblock {\em submitted}.

\bibitem{popovic07b}
M.~Krupa, N.~Popovi\'c, H.G. Rotstein, and N.~Kopell.
\newblock Mixed-mode oscillations in a three time-scale model for the
  dopaminergic neuron.
\newblock {\em Chaos}, 2007.

\bibitem{Acker03}
C.D. Acker, N.~Kopell, and J.A. White.
\newblock Synchronization of strongly coupled excitatory neurons: Relating
  network behavior to biophysics.
\newblock {\em J. Comput. Neurosci.}, 15:71 -- 90, 2003.

\bibitem{Alonso87}
A.~Alonso and E.~Garcia-Austt.
\newblock Neuronal sources of the theta rhythm in the entorhinal cortex of the
  rat. i. laminar distribution of theta filed potentials.
\newblock {\em Expl Brain Res.}, 67:493--501, 1987.

\bibitem{Mitchell80}
S.J. Mitchell and Ranck J.~B. Jr.
\newblock Generation of the theta rhythm in medial entorhinal cortex of freely
  moving rats.
\newblock {\em Brain Res.}, 189:49 --66, 1980.

\bibitem{kn:win1}
J.~Winson.
\newblock Loss of hippocampal theta rhythm results in spatial memory deficit in
  the rat.
\newblock {\em Science}, 201:160--163, 1978.

\bibitem{kn:buz2}
G.~Buzs\'aki.
\newblock Two-stage model of memory trace formation: a role for 'noisy' brain
  states.
\newblock {\em Neuroscience}, 31:551--570, 1989.

\bibitem{Schmitz98}
D.~Schmitz, T.~Gloveli, J.~Behr, T.~Dugladze, and U.~Heinemann.
\newblock Subthreshold membrane potential oscillations in neurons of deep
  layers of the entorhinal cortex.
\newblock {\em Neuroscience}, 85(4):999--1004, 1998.

\bibitem{Agrawal01}
N.~Agrawal, B.N. Hamam, J.~Magistretti, A.~Alonso, and D.S. Ragsdale.
\newblock Persistent sodium channel activity mediates subthreshold membrane
  potential oscillations and low-threshold spikes in rat entorhinal cortex
  layer v neurons.
\newblock {\em Neuroscience}, 102(1):53--64, 2001.

\bibitem{kn:dicalo3}
C.~T. Dickson, J.~Magistretti, M.~H Shalinsky, E.~Frans\'{e}n, M.~Hasselmo, and
  A.~A. Alonso.
\newblock Properties and role of \( i_{h} \) in the pacing of subthreshold
  oscillation in entorhinal cortex layer ii neurons.
\newblock {\em J. Neurophysiol.}, 83:2562--2579, 2000.

\bibitem{kn:alokli2}
R.~M. Klink and A~Alonso.
\newblock Ionic mechanisms for the subthreshold oscillations and differential
  electroresponsiveness of medial entorhinal cortex layer ii neurons.
\newblock {\em J. Neurophysiol.}, 70:128--143, 1993.

\bibitem{kn:clekop1}
R.~Clewley, H.~G. Rotstein, and N.~Kopell.
\newblock A computational tool for the reduction of nonlinear ode systems
  possessing multiple scales.
\newblock {\em Multiscale Model. Simul.}, 4:732--759, 2005.

\bibitem{LS}
R.~Larter and C.G. Steinmetz.
\newblock Chaos via mixed mode oscillations.
\newblock {\em Phil. Trans. R. Soc. London}, A337:241--298, 1991.

\bibitem{arne}
A.~Arn\'eado, F.~Argoul, J.~Elezgaray, and P.~Richetti.
\newblock Homoclinic chaos in chemical systems.
\newblock {\em Physica D}, 62:134--168, 1993.

\bibitem{Guckenheimer1}
J.~Guckenheimer, R.~Harris-Warrick, A.~Peck, and A.~Willms.
\newblock Bifurcation, bursting, and spike frequency adaptation.
\newblock {\em Journal of Computational Neuroscience}, 4:257--277, 1997.

\bibitem{Guckenheimer2}
J.~Guckenheimer and A.~Willms.
\newblock Asymptotic analysis of subcritical Hopf-homoclinic bifurcation.
\newblock {\em Physica D}, 139:195--216, 2000.

\bibitem{Benoit90}
E~Benoit.
\newblock Canards et enlacements.
\newblock {\em Publ. Inst. Hautes Etudes Sci.}, 72:63--91, 1990.

\bibitem{Szmolyan01}
P.~Szmolyan and M.~Wechselberger.
\newblock Canards in $\mathbb{R}^3$.
\newblock {\em JDE}, 177(2):419--453, 2001.

\bibitem{Wechselberger05}
M.~Wechselberger.
\newblock Existence and bifurcation of canards in $\mathbb{R}^3$ in the case of
  a folded node.
\newblock {\em SIAM J. Appl. Dyn. Sys.}, 4:101--139, 2005.

\bibitem{BKW06}
M.~Brons, M.~Krupa, and M.~Wechselberger.
\newblock Mixed mode oscillations due to the generalized canard phenomenon.
\newblock {\em Fields Institute Communications}, 49:39--63, 2006.

\bibitem{popovic07a}
M.~Krupa, N.~Popovi\'c, and N.~Kopell.
\newblock Mixed-mode oscillations in three time-scale systems: a prototypical
  example.
\newblock {\em preprint}, 2007.

\bibitem{Krupa01}
M.~Krupa and P.~Szmolyan.
\newblock Extending geometric singular perturbation theory to non-hyperbolic
  points--fold and canard points in two dimensions.
\newblock {\em SIAM J. Math. Anal.}, 33(2):286--314, 2001.

\bibitem{Dumortier96}
F.~Dumortier and R.~Roussarie.
\newblock Canard cycles and center manifolds.
\newblock {\em Memoirs of the AMS}, 121 (577), 1996.

\bibitem{Eckhaus83}
W.~Eckhaus.
\newblock Relaxation oscillations including a standard chase on french ducks,
  in {\em asymptotic analysis ii}.
\newblock {\em Springer Lecture Notes in Mathematics}, 985:449--494, 1983.

\bibitem{Hodgkin52}
A.L. Hodgkin and A.F. Huxley.
\newblock A quantitative description of membrane current and its application to
  conduction and excitation in a nerve.
\newblock {\em Journal of Physiology}, 117:165--181, 1952.

\bibitem{kn:izh2}
E.~Izhikevich.
\newblock {\em Dynamical Systems in Neuroscience: The geometry of excitability
  and bursting}.
\newblock MIT Press (Cambridge, Massachusetts), 2006.

\bibitem{Richardson03}
M.J.E. Richardson, N.~Brunel, and V.Hakim.
\newblock From subthreshold to firing rate resonance.
\newblock {\em J. Neurophysiol}, 89:2538--2554, 2003.

\bibitem{Shilnikov65}
L.P. Shilnikov.
\newblock A case of the existence of a denumerable set of periodic motions.
\newblock {\em Sov. Math. Dokl.}, 6, 1965.

\bibitem{Fenichel71}
N.~Fenichel.
\newblock Persistence and smoothness of invariant manifolds for flows.
\newblock {\em Indiana Univ. Mathematics Journal}, 21:193--226, 1971.

\bibitem{Benoit81}
E.~Benoit, J.~F. Callot, F~Dienner, and M.~Dienner.
\newblock Chasse au canard.
\newblock {\em Collectanea Mathematica}, 31-32 (1-3):37--119, 1981.

\bibitem{Diener94}
M.~Diener.
\newblock Regularizing microscopes and rivers.
\newblock {\em SIAM J. Appl. Math.}, 25:148--173, 1994.

\bibitem{Carr81}
J.~Carr.
\newblock {\em Applications of centre manifold theory}.
\newblock Springer-Verlag, 1981.

\bibitem{Chrobak94}
J.~J. Chrobak and G.~Buzs\'aki.
\newblock Selective activation of deep layer (v-vi) retrohippocampal cortical
  neurons during hippocampal sharp waves in the behaving rat.
\newblock {\em J. Neurosci.}, 14:6160--6170, 1994.

\bibitem{Jones90}
R.S.G. Jones and J.D.C. Lambert.
\newblock Synchronous discharges in the rat entorhinal cortex \emph{in vitro}:
  site of initiation and the role of excitatory amino acid receptors.
\newblock {\em Neuroscience}, 34:657--670, 1990.

\bibitem{Egorov02}
A.V. Egorov, B.N. Hamam, E.~Fransen, M.E. Hasselmo, and A.A. Alonso.
\newblock Graded persistent activity in entorhinal cortex neurons.
\newblock {\em Nature}, 420:173--178, 2002.

\bibitem{Koper92}
M.T.M. Koper and P.~Gaspard.
\newblock The modeling of mixed-mode and chaotic oscillations in
  electrochemical systems.
\newblock {\em J. Chem. Phys}, 96:7797--7813, 1992.

\bibitem{Fenichel79}
N.~Fenichel.
\newblock Geometric singular perturbation theory.
\newblock {\em J. Diff. Eq.}, 31:53--98, 1979.

\bibitem{Jones95}
C.K.R.T Jones.
\newblock {\em Geometric singular perturbation theory, in \emph{Dynamical
  Systems}}, volume 1609.
\newblock Springer Lecture Notes in Mathematics, 1995.

\bibitem{Rajesh00}
Rajesh S. and G.~Ananthakrishna.
\newblock Incomplete approach to homoclinicity in a model with bent-slow
  geometry.
\newblock {\em Physica D}, 140:193--212, 2000.

\end{thebibliography}

\end{document}